\documentclass[11pt,letterpaper,reqno]{amsart}

\usepackage[T1]{fontenc}
\usepackage[letterpaper,margin=1in]{geometry}
\usepackage{amsmath,amssymb,amsthm,mathtools}
\numberwithin{equation}{section}
\usepackage{booktabs,tabularx,array,float}
\usepackage{microtype}
\usepackage[ruled,vlined,linesnumbered]{algorithm2e}
\usepackage[hidelinks]{hyperref}
\hypersetup{
  pdftitle={Optimal Rational Spherical Designs on the Circle},
  pdfauthor={Zhuo Cheng, Yeyao Hu, Deyu Yu},
  pdfsubject={Optimal rational circle designs, universal local divisibility in arbitrary strength, exact sigma and orbit trades, and nonuniqueness},
  pdfkeywords={Rational spherical designs, Gaussian integers, optimal cardinality, torus immersions}
}

\newtheorem{theorem}{Theorem}[section]
\newtheorem{proposition}[theorem]{Proposition}
\newtheorem{lemma}[theorem]{Lemma}
\newtheorem{corollary}[theorem]{Corollary}

\theoremstyle{definition}

\theoremstyle{remark}
\newtheorem{remark}[theorem]{Remark}
\newtheorem{problem}[theorem]{Problem}
\newcommand{\ii}{\mathrm{i}}
\newcommand{\Z}{\mathbb{Z}}
\newcommand{\Q}{\mathbb{Q}}
\newcommand{\C}{\mathbb{C}}
\newcommand{\T}{\mathbb{T}}

\newcommand{\Norm}{\operatorname{N}}
\newcommand{\RePart}{\operatorname{Re}}
\newcommand{\ImPart}{\operatorname{Im}}
\newcommand{\Korb}{\mathcal{K}}
\newcommand{\Orb}{\mathcal{O}}
\newcommand{\Efour}{\mathcal{E}_4}
\newcommand{\Efive}{\mathcal{E}_5}
\newcommand{\Hfour}{H_4}
\newcommand{\Hfive}{H_5}
\newcommand{\muFour}{\mu_4}
\newcommand{\NQ}{N_{\Q}}

\SetAlgoCaptionSeparator{:}
\SetAlCapSkip{0.4em}
\SetAlCapHSkip{0.15em}
\SetAlgoNlRelativeSize{-1}
\SetNlSty{textbf}{}{}
\SetInd{0.55em}{0.85em}
\DontPrintSemicolon
\SetKwInput{Input}{Input}
\SetKwInput{Output}{Output}
\SetKw{Continue}{continue}
\SetKw{End}{end}
\SetKw{Return}{return}
\SetKw{Verify}{verify}
\SetKw{Emit}{emit}
\SetKwBlock{Begin}{begin}{end}
\SetCommentSty{itshape}

\newcommand{\Audit}{\mathsf{Audit}}

\newcommand{\UNKNOWN}{\textnormal{\textsc{Unknown}}}

\title{Optimal Rational Spherical Designs on the Circle}

\author{Zhuo Cheng}
\address{School of Mathematics and Statistics, HNP-LAMA, Central South
  University, Changsha, Hunan 410083, P. R. China}
\email{zhuo@outlook.cz}

\author{Yeyao Hu}
\address{School of Mathematics and Statistics, HNP-LAMA, Central South
  University, Changsha, Hunan 410083, P. R. China}
\email{huyeyao@gmail.com}

\author{Deyu Yu}
\address{School of Mathematics and Statistics, HNP-LAMA, Central South
  University, Changsha, Hunan 410083, P. R. China}
\email{242107002@csu.edu.cn}

\date{November 2026}

\begin{document}

\begin{abstract}
The minimum numbers of points in rational equal-weight spherical designs
on the circle are $2304$, $4608$, $9216$, and $9216$ at strengths $4$,
$5$, $6$, and $7$, respectively.  We construct two nonisometric optimal
designs at each strength. At arbitrary strength $t$, we obtain an explicit cardinality
divisor $D_t$ with $\log D_t=(1/2+o(1))t\log t$.
Applied to the antipodal pairs of a $4608$-point design, Gromov's torus
construction gives a flat minimal immersion of $\T^2$ into $S^{4607}$
with $2304$ distinct integral frequencies, attaining Petrunin's sharp
Euclidean normal-curvature bound. 
\end{abstract}

\maketitle

\section{Introduction and main results}\label{sec:introduction}

Let $\sigma_{d-1}$ denote the rotation-invariant probability measure on
$S^{d-1}\subset\mathbb{R}^d$, where $d\geq2$.  A nonempty finite set
$X\subset S^{d-1}$ is a spherical $t$-design if
\[
 \frac1{|X|}\sum_{p\in X}f(p)
 =\int_{S^{d-1}}f\,d\sigma_{d-1}
\]
for every real polynomial $f$ of total degree at most $t$.  Introduced
by Delsarte, Goethals, and Seidel~\cite{DGS}, spherical designs are
equal-weight cubature formulas on the sphere.  A design is rational when
$X\subset S^{d-1}\cap\Q^d$.  We study rational designs on $S^1$, always
with distinct nodes and equal weights.  Classical results on circle
$t$-designs, without a rationality restriction, include Hong's
classification at small cardinalities~\cite{Hong82}.  The general real
existence theorem of Seymour and Zaslavsky~\cite{SZ84} also places no
arithmetic restriction on the nodes.

Requiring rational coordinates is different from asking for an explicit
construction.  Cui, Xia, and Xiang~\cite{CXX} distinguish these questions
and prove rational existence on suitable open connected spaces, as well
as spherical existence over real abelian extensions of $\Q$.
They identify the existence of a rational
spherical $4$-design on $S^1$ as the smallest open case.
Xiang's explicit constructions over real abelian
extensions~\cite{Xiang22} and Misawa's explicit $5$- and $7$-designs
\cite{Mis26} likewise do not settle this rational distinct-node problem.

Gromov's construction of curvature-controlled tori already appears in
\cite[Sections~2--3]{Gro22}.  In his later lectures
\cite[Section~13]{Gro25}, he explains the role of rational fourth-moment
identities and contrasts the regular pentagon with the lack of a simple
rational circle design.  The geometric construction and its optimality,
proved by Petrunin~\cite{Pet24}, are prior results.  Here the arithmetic
problem is to give distinct rational nodes and determine their minimum
cardinality at four consecutive strengths.

Rational moment equations also lead to Diophantine problems.  In their
study of Hilbert--Kamke equations for classical orthogonal-polynomial
measures, Mishima, Lu, Sawa, and Uchida~\cite{MLSU26} use rational
equal-weight designs for the Chebyshev measure
$dx/(\pi\sqrt{1-x^2})$ on $(-1,1)$ to obtain solutions of the
Prouhet--Tarry--Escott problem.  Matsumura and Sawa~\cite{MS25} construct
parametric rational ideal solutions of its two-dimensional version from
pairs of six-point rational $5$-designs on the ellipse $x^2+xy+y^2=1$.
For rational circle designs, clearing denominators places the nodes on
an integer shell $a^2+b^2=H^2$.  One cannot simply take the whole shell.
Bannai and Miezaki~\cite{BM10} proved that no nonempty positive-norm
shell of $\Z^2$ is a spherical $4$-design after normalization;
Bannai, Miezaki, and Yudin~\cite{BMY} subsequently gave elementary
proofs of this obstruction.  We instead select proper subsets of the
shells.

For every integer $t\geq1$, define
\[
 \NQ(t)=\min\{|Y|:Y\subset S^1\cap\Q^2
 \text{ is an equal-weight spherical $t$-design}\},
\]
with the convention $\NQ(t)=+\infty$ if no such finite design exists.

Throughout, ``optimal'' means having the least cardinality among rational
equal-weight designs.  A design has exact strength $t$ if it is a
$t$-design but not a $(t+1)$-design.  Identifying $S^1$ with the unit
circle in $\C$, we write $P_k(X)=\sum_{z\in X}z^k$.
Prescribing the entire set of degrees for which $P_k$ vanishes is a
different problem, studied for real circle configurations by Misawa
and Nishimura~\cite{MN26}; no rationality condition is imposed there.
Our use of exact strength concerns only the initial vanishing segment
and the next moment.

Our first $4$-design is obtained by deleting conjugate pairs from a
Gaussian shell.  Some $D_4$-orbits are kept in full; from each remaining
orbit we remove one conjugate pair. 
The resulting set is conjugation-invariant but not antipodal. We also give a
second design with a different fifth moment.

\begin{theorem}\label{thm:optimal4}
There exists a set
\[
 X_4\subset S^1\cap\Q^2
\]
of $2304$ distinct points which is an equal-weight spherical $4$-design but
not a spherical $5$-design.  Its strength is exactly $4$, and
\[
 X_4\subset \Hfour^{-1}\Z^2\cap S^1,
 \qquad \Hfour=48\,612\,265.
\]
\end{theorem}

\begin{theorem}
\label{thm:second4-intro}
There exists a set
\[
 \widetilde X_4\subset \Hfour^{-1}\Z^2\cap S^1
\]
of $2304$ distinct points and exact strength $4$.  It is obtained from the
selector defining $X_4$ by an exact $157$-packet sign trade.  No isometry of
$S^1$ maps $X_4$ onto $\widetilde X_4$.
\end{theorem}

The lower bound uses only $P_4=0$: any list of $N$ rational points on
$S^1$ with this property satisfies $2304\mid N$, even when points are
repeated and no symmetry is imposed.  Both constructions attain the bound.

\begin{corollary}\label{cor:NQ}
One has
\[
 \NQ(4)=2304.
\]
Consequently both $X_4$ and $\widetilde X_4$ are globally
cardinality-optimal rational equal-weight spherical $4$-designs.
\end{corollary}

For strengths $5$, $6$, and $7$, we work on a larger Gaussian shell and give a second design.

\begin{theorem}\label{thm:main}
There exists a set
\[
 X_5\subset S^1\cap\Q^2
\]
of $4608$ distinct points which is an equal-weight spherical $5$-design but
not a spherical $6$-design.  Its strength is exactly $5$.  Moreover,
\[
 X_5\subset \Hfive^{-1}\Z^2\cap S^1,
 \qquad \Hfive=12\,882\,250\,225,
\]
and $X_5$ is invariant under $z\mapsto-z$ and
$z\mapsto\overline z$.
\end{theorem}

The even moments of $X_5$ can also be used in Gromov's torus construction
\cite{Gro22,Gro25}, as formulated in
Chodosh--Li~\cite[Remark~11]{CL}.  Taking one representative from each
antipodal pair gives an equal-amplitude immersion with distinct integral
frequencies on a single shell.  Its curvature attains Petrunin's lower
bound~\cite{Pet24}.  The contribution here is the explicit certified
realization with $2304$ distinct frequencies on the specified integral
shell, not the curvature value or the general design-to-torus mechanism;
see Theorem~\ref{thm:torus-intro}.

To prove optimality at strength $5$, we lift a fourth-order congruence
using four parameters.  This gives $4608\mid N$ without symmetry or
distinctness assumptions.

\begin{corollary}\label{cor:NQ5}
One has
\[
 \NQ(5)=4608.
\]
Thus $X_5$ is globally cardinality-optimal among rational equal-weight
spherical $5$-designs on the circle.
\end{corollary}

Let $D_4$ denote the order-eight group generated by a quarter-turn and complex
conjugation.

\begin{theorem}\label{thm:parent-intro}
There exists a $D_4$-invariant rational spherical $7$-design
\[
 X_7=X_5\sqcup\ii X_5
\]
of $9216$ distinct points.  Its exact strength is $7$, and the same point
set is globally cardinality-optimal among rational equal-weight spherical
$6$-designs and among rational equal-weight spherical $7$-designs.  Moreover,
every $D_4$-invariant rational spherical $4$-design has cardinality
divisible by $9216$, so within that symmetry class the same lower bound
already follows from the fourth moment.
\end{theorem}

At strength $6$, we gives $9216\mid N$ under the same
general assumptions.  The bound also holds at strength $7$.
For arbitrary strength, Theorem~\ref{thm:universal-Dt} gives an explicit
divisor $D_t$ of every rational equal-weight $t$-design cardinality,
with $\log D_t=(1/2+o(1))t\log t$ as $t\to\infty$.

\begin{corollary}
\label{cor:NQ67}
One has
\[
 \NQ(6)=\NQ(7)=9216.
\]
Thus the exact-strength-$7$ set $X_7$ is globally cardinality-optimal both as
a rational equal-weight spherical $6$-design and as a rational equal-weight
spherical $7$-design.
\end{corollary}

\begin{theorem}
\label{thm:nonunique-intro}
There are sets
\[
 \widetilde{X}_5\subset \Hfive^{-1}\Z^2\cap S^1,
 \qquad
 \widetilde{X}_7=\widetilde{X}_5\sqcup\ii\widetilde{X}_5,
\]
with $|\widetilde{X}_5|=4608$ and $|\widetilde{X}_7|=9216$, of exact
strengths $5$ and $7$, respectively.  No isometry of $S^1$ maps $X_5$ onto
$\widetilde{X}_5$, and no isometry maps $X_7$ onto $\widetilde{X}_7$.
Consequently globally cardinality-optimal rational equal-weight spherical
designs at strengths $5$, $6$, and $7$ are not unique, even up to rotations
and reflections.
\end{theorem}

The paper is organized as follows.
Section~\ref{sec:fourier-divisibility} proves
Theorems~\ref{thm:universal-divisibility},
\ref{thm:fifth-divisibility},
\ref{thm:sixth-divisibility}, and
\ref{thm:universal-Dt}.
Using Fourier moments and local arithmetic congruences,
we establish cardinality divisors at strengths
$4$, $5$, $6$, and $7$, together with a universal divisor
$D_t$ for arbitrary strength.

Section~\ref{sec:optimal567} proves
Theorems~\ref{thm:main},
\ref{thm:parent-intro}, and
\ref{thm:nonunique-intro}.
We construct optimal rational designs at strengths
$5$, $6$, and $7$ from a common Gaussian integer shell
and obtain nonisometric optimal configurations
through exact orbit trades.

Section~\ref{sec:optimal2304} proves
Theorems~\ref{thm:optimal4} and
\ref{thm:second4-intro}.
We construct two nonisometric optimal rational
$4$-designs by conjugate-pair deletions and an exact
sign trade, and obtain a family of exact-strength-$4$
designs with $2304r$ points for every $r\geq1$.

Section~\ref{sec:applications} proves
Theorem~\ref{thm:torus-intro} by applying Gromov's
torus construction to an optimal rational $5$-design.
It also discusses the possible infinitude of optimal
designs modulo $O(2)$.

Appendix~\ref{app:proof-search} describes five
computational searches and their exact acceptance tests.
Appendix~\ref{app:certificates} contains the
construction certificates and their independent verification.

\section{Fourier moments and arithmetic divisibility}
\label{sec:fourier-divisibility}

On the circle, the design equations are equivalent to the vanishing of
finitely many Fourier moments.  Combining these equations with local
arithmetic in $\Q(\ii)$ gives the divisors $2304$, $4608$, and $9216$.
We first prove the low-degree bounds, then discuss the additional
restrictions imposed by symmetry and the general divisor $D_t$.

\subsection{Arithmetic notation}\label{subsec:arithmetic-notation}

We use standard facts about Gaussian factorization, quadratic splitting,
and finite-field norms; background can be found in Ireland and Rosen
\cite[Chapters~1--2, 7, and~13]{IR90}.  The specific cardinality
congruences are proved below.

For a rational prime $p$, write $v_p$ for the valuation on $\Q$ normalized
by $v_p(p)=1$, and set $v_p(0)=+\infty$.  Thus
\[
 \Z_{(p)}=\{x\in\Q:v_p(x)\geq0\},\qquad
 \Z_{(p)}^\times=\{x\in\Q:v_p(x)=0\}.
\]
We interpret $a\equiv b\pmod{p^m}$ in $\Z_{(p)}$ as
$v_p(a-b)\geq m$.  For odd $n\geq1$, the symbol
$\bigl(\frac an\bigr)$ denotes the Jacobi symbol, and for $n=p$ prime
it is the Legendre symbol.  We use the quadratic Dirichlet character
of conductor $4$,
\[
 \chi_{-4}(n)=
 \begin{cases}
  \displaystyle\left(\frac{-1}{n}\right)=(-1)^{(n-1)/2},&2\nmid n,\\[2pt]
  0,&2\mid n.
 \end{cases}
\]
An odd prime $p$ splits in $\Q(\ii)$ when $\chi_{-4}(p)=1$ and is
inert when $\chi_{-4}(p)=-1$.  Write
\[
 \mathcal P_-=
 \{p\text{ prime}:\chi_{-4}(p)=-1\}
 =\{p\text{ prime}:p\equiv3\pmod4\}.
\]

\subsection{Fourier moments and symmetry}\label{subsec:fourier-symmetry}

Identify $(x,y)\in S^1$ with $z=x+\ii y\in\C$, and write
\[
 P_k(X)=\sum_{z\in X}z^k.
\]

The following is the circle form of the standard harmonic criterion for
spherical designs~\cite{DGS}; see also Hong~\cite[Lemma~1]{Hong82}.
We include the short proof to fix the normalization.

\begin{lemma}\label{lem:fourier}
A finite set $X\subset S^1$ is an equal-weight spherical $t$-design if and
only if
\[
 P_k(X)=0\qquad(1\leq k\leq t).
\]
\end{lemma}

\begin{proof}
On $S^1$,
\[
 x=\frac{z+z^{-1}}2,
 \qquad
 y=\frac{z-z^{-1}}{2\ii}.
\]
Thus a polynomial of degree at most $t$ restricts to a Laurent polynomial
$\sum_{k=-t}^t c_kz^k$.  For $k\neq0$, the circular average of $z^k$
is zero.  The negative moments are conjugates of the positive moments,
so the stated conditions imply the design identity.  Conversely, apply
that identity to the real and imaginary parts of $(x+\ii y)^k$.
\end{proof}

For the invariant-cubature principle behind the symmetry reduction,
see Sobolev~\cite{Sob62} and the finite-orthogonal-group formulation in
\cite[Theorem~3.2]{MS25}.  In the groups used here, the reduction follows
by explicit orbit sums.
Let
\[
 K=\langle z\mapsto-z,\ z\mapsto\overline z\rangle
 \cong C_2\times C_2.
\]
For $w\in S^1$, define
\[
 \Korb(w)=\{w,-w,\overline w,-\overline w\},
\]
and, with $\muFour=\Z[\ii]^\times=\{1,-1,\ii,-\ii\}$,
\[
 \Orb(w)=\{uw,u\overline w:u\in\muFour\}
 =\Korb(w)\cup\Korb(\ii w).
\]

\begin{lemma}\label{lem:orbit-moments}
If $w=a+\ii b\in S^1$, $ab\neq0$, and $|a|\neq|b|$, then
$|\Korb(w)|=4$, $|\Orb(w)|=8$, and
\[
 \sum_{z\in\Korb(w)}z^k=
 \begin{cases}
 0,&k\text{ odd},\\
 4\RePart(w^k),&k\text{ even},
 \end{cases}
\]
while
\[
 \sum_{z\in\Orb(w)}z^k=
 \begin{cases}
 0,&4\nmid k,\\
 8\RePart(w^k),&4\mid k.
 \end{cases}
\]
If $w=(A+\ii B)/C$ and $A^2+B^2=C^2$, then
\begin{equation}\label{eq:Rew4}
 \RePart(w^4)=1-\frac{8A^2B^2}{C^4}.
\end{equation}
\end{lemma}

\begin{proof}
The assumptions on $a,b$ rule out the axes and diagonals, so neither
orbit has fewer points.  The first moment formula follows by summing
over antipodal pairs and their conjugates.  For the second, use
\[
 \sum_{u^4=1}(uw)^k=w^k\sum_{u^4=1}u^k.
\]
Finally,
\[
 \RePart((A+\ii B)^4)
 =A^4-6A^2B^2+B^4
 =(A^2+B^2)^2-8A^2B^2.
\]
\end{proof}

\begin{corollary}\label{cor:automatic}
A $K$-invariant spherical $4$-design is automatically a spherical
$5$-design.  A $D_4$-invariant spherical $4$-design is automatically a
spherical $7$-design.
\end{corollary}

\begin{proof}
Antipodality gives vanishing of every odd moment.  Under a quarter-turn,
$P_k=\ii^kP_k$, so the moment also vanishes whenever its degree is not
divisible by $4$.
\end{proof}

\subsection{The fourth-moment divisor}\label{subsec:fourth-divisibility}

\begin{lemma}\label{lem:common-denominator}
Let $z_1,\ldots,z_N\in S^1\cap\Q(\ii)$.  There exists a positive integer $C$, coprime to $6$, such that
\[
 z_j=\frac{A_j+\ii B_j}{C},
 \qquad A_j,B_j\in\Z,
 \qquad A_j^2+B_j^2=C^2
\]
for every $j$.  Moreover,
\[
 12\mid A_jB_j
\]
for every $j$.
\end{lemma}

\begin{proof}
Choose primitive representations
\[
 z_j=\frac{a_j+\ii b_j}{c_j},\qquad
 a_j^2+b_j^2=c_j^2,\qquad
 \gcd(a_j,b_j,c_j)=1,\qquad c_j>0.
\]
If $2\mid c_j$, reduction modulo $4$ forces $2\mid a_j,b_j$, a
contradiction.  Now let $p$ be an odd prime divisor of $c_j$.
Primitivity gives $p\nmid a_jb_j$, and hence
$a_j^2\equiv-b_j^2\pmod p$ yields
\[
 \left(\frac{-1}{p}\right)
 =\left(\frac{(a_jb_j^{-1})^2}{p}\right)=1,
\]
where $b_j^{-1}$ is taken modulo $p$.  Thus every prime divisor of $c_j$
is split.  Set $C=\operatorname{lcm}(c_1,\ldots,c_N)$; then
\[
 v_2(C)=0,\qquad v_p(C)=0\quad(p\in\mathcal P_-).
\]
In particular, $\gcd(C,6)=1$ since $\chi_{-4}(3)=-1$.
Set $A_j=(C/c_j)a_j$ and $B_j=(C/c_j)b_j$.

Since $C$ is odd, reduction of $A_j^2+B_j^2=C^2$ modulo $2$ gives
$\min\{v_2(A_j),v_2(B_j)\}=0$.  Reduction modulo $8$ then shows that
$v_2(A_jB_j)\geq2$.  Modulo $3$, exactly one coordinate vanishes,
because $C$ is a unit and the only nonzero square is $1$.  Thus
$v_3(A_jB_j)\geq1$, proving $2^2\cdot3\mid A_jB_j$.
\end{proof}

In particular, every rational circle point $z=x+\ii y$ satisfies
\begin{equation}\label{eq:local-coordinate-valuations}
 \begin{gathered}
 x,y\in\Z_{(p)}\quad(p=2\text{ or }p\in\mathcal P_-),\\
 \min\{v_2(x),v_2(y)\}=0,\qquad v_2(xy)\geq2,\qquad
 v_3(xy)\geq1.
 \end{gathered}
\end{equation}
In the $2$-adic arguments below, we distinguish the two possibilities
$(x,y)\in\Z_{(2)}^\times\times4\Z_{(2)}$ and
$(x,y)\in4\Z_{(2)}\times\Z_{(2)}^\times$.

\begin{theorem}\label{thm:universal-divisibility}
Let $z_1,\ldots,z_N\in S^1\cap\Q(\ii)$, with repetitions allowed.  If
\[
 \sum_{j=1}^N z_j^4=0,
\]
then
\[
 2304\mid N.
\]
In particular, every rational equal-weight spherical $4$-design on $S^1$ has
cardinality divisible by $2304$.
\end{theorem}

\begin{proof}[Proof of Theorem~\ref{thm:universal-divisibility}]
Choose $C$ as in Lemma~\ref{lem:common-denominator} and write
\[
 z_j=\frac{A_j+\ii B_j}{C},
 \qquad
 q_j=\frac{A_jB_j}{12}\in\Z,
 \qquad
 \Delta_j=A_j^2-B_j^2.
\]
The coordinate conditions imply $\Delta_j\in\Z_{(2)}^\times$.
Expanding the fourth power gives
\begin{align}
 \RePart(A_j+\ii B_j)^4&=C^4-1152q_j^2,\label{eq:fourth-real-arithmetic}\\
 \ImPart(A_j+\ii B_j)^4&=48q_j\Delta_j.\label{eq:fourth-imag-arithmetic}
\end{align}
The real part of $\sum_j z_j^4=0$ is
\[
 NC^4=1152\sum_{j=1}^N q_j^2.
\]
As $v_2(C)=v_3(C)=0$ and $1152=2^7 3^2$, this gives
$v_2(N)\geq7$ and $v_3(N)\geq2$.  Writing $N=1152m$, we obtain
\[
 \sum_{j=1}^N q_j^2=mC^4.
\]
Modulo $2$, we have $q_j^2\equiv q_j$ and $C^4\equiv1$, whence
\begin{equation}\label{eq:q-parity-m}
 \sum_{j=1}^N q_j\equiv m\pmod2.
\end{equation}
The imaginary part gives, by \eqref{eq:fourth-imag-arithmetic},
\[
 \sum_{j=1}^N q_j\Delta_j=0.
\]
Since $\Delta_j\equiv1\pmod2$ for every $j$,
\[
 \sum_{j=1}^N q_j\equiv0\pmod2.
\]
Comparison with \eqref{eq:q-parity-m} gives $v_2(m)\geq1$.  Hence
\[
 v_2(N)\geq8,\qquad v_3(N)\geq2,
\]
which is equivalent to $2^8 3^2=2304\mid N$.
\end{proof}

\subsection{Fifth-design divisibility and four-parameter lifting}
\label{subsec:fifth-divisibility}

\begin{theorem}\label{thm:fifth-divisibility}
Let $z_1,\ldots,z_N\in S^1\cap\Q(\ii)$, with repetitions allowed.  If
\[
 \sum_{j=1}^N z_j^k=0\qquad(k=1,3,4,5),
\]
then
\[
 4608\mid N.
\]
In particular, every rational equal-weight spherical $5$-design on $S^1$ has
cardinality divisible by $4608$.
\end{theorem}

To gain the extra factor of $2$ at strength $5$, we add odd terms of
degree at most $5$ to a fourth-order congruence.  Let
$z=x+\ii y\in S^1\cap\Q(\ii)$.  By \eqref{eq:local-coordinate-valuations},
we may write $xy=4u$ with $u\in\Z_{(2)}$ and
$\Delta=x^2-y^2\in\Z_{(2)}^\times$.  Set
\[
 \Phi_4(z)=2\RePart(z^4)+16\ImPart(z^4).
\]
The identities $\RePart(z^4)-1=-8x^2y^2$ and
$\ImPart(z^4)=4xy\Delta$ give
\begin{equation}\label{eq:phi4-seed}
 \Phi_4(z)-2=2^8u(\Delta-u)\in2^9\Z_{(2)}.
\end{equation}
The divisibility follows from $\Delta\equiv1\pmod2$, which gives
$v_2(u(\Delta-u))\geq1$.  An odd polynomial of degree at most $5$ has
zero average on the circle and on every $5$-design.  We seek such a
polynomial to cancel the remaining residue in
$2^9\Z_{(2)}/2^{10}\Z_{(2)}\cong\mathbb F_2$.  There is no parity
obstruction to this step: each residue class is its own negative.

For integers $a,b,c,d$, define, on the unit circle,
\begin{equation}\label{eq:phi-four-param-family}
 \Phi_{a,b,c,d}(x+\ii y)
 =\Phi_4(x+\ii y)+16\bigl(axy^2+bxy^4+cx^2y+dx^2y^3\bigr).
\end{equation}
All four correction terms are odd, have degree at most $5$, and vanish
at the axial points.  At those points, $\Phi_4=2$.

\begin{lemma}\label{lem:four-param-lift}
Suppose that $a,b,c,d\in\Z$ satisfy
\begin{equation}\label{eq:phi-parameter-conditions}
 \begin{gathered}
 a\equiv1\pmod2,\qquad c+4a\equiv0\pmod{16},\qquad
 a+b\equiv4\pmod8,\\
 c+d\equiv\frac{a+b}{4}\pmod4.
 \end{gathered}
\end{equation}
Then $v_2(\Phi_{a,b,c,d}(z)-2)\geq10$ for every
$z\in S^1\cap\Q(\ii)$; equivalently,
$\Phi_{a,b,c,d}(z)\equiv2\pmod{2^{10}}$.
\end{lemma}

\begin{proof}
On $x^2+y^2=1$, expansion gives
\begin{equation}\label{eq:phi-parameter-factor}
 \Phi_{a,b,c,d}(x+\ii y)-2
 =16xy\bigl(4(x^2-y^2)-xy+ay+by^3+cx+dxy^2\bigr).
\end{equation}
Write $h=c/4$ and $q=(a+b)/4$.  The parameter conditions imply
$h,q\in\Z\cap\Z_{(2)}^\times$, together with $h\equiv-a\pmod4$ and
$c+d\equiv q\pmod4$.

First suppose $v_2(y)\geq2$ and write $y=4u$.
Then $x\in\Z_{(2)}^\times$ and $x^2=1-16u^2$.
Substitution in \eqref{eq:phi-parameter-factor} gives
\begin{equation}\label{eq:phi-parameter-y-case}
 \Phi_{a,b,c,d}(x+4\ii u)-2
 \equiv256xu\bigl(1+hx+u(a-x)\bigr)\pmod{1024}.
\end{equation}
Modulo $4$, the factor after $256x$ is
\[
 u\bigl(1+hx+u(a-x)\bigr)\equiv
 \begin{cases}
 (a-1)u(u-1),&x\equiv1\pmod4,\\
 (a+1)u(u+1),&x\equiv-1\pmod4.
 \end{cases}
\]
Since $v_2(a\pm1)\geq1$ and $v_2(u(u\pm1))\geq1$, both expressions
have valuation at least $2$.  Including the factor $256=2^8$ gives
$v_2(\Phi_{a,b,c,d}-2)\geq10$.

Now suppose $v_2(x)\geq2$.  Write $x=4u$; then
$y\in\Z_{(2)}^\times$ and $y^2=1-16u^2$.  The same substitution yields
\begin{equation}\label{eq:phi-parameter-x-case}
 \Phi_{a,b,c,d}(4u+\ii y)-2
 \equiv256uy\bigl(-1+qy+u(c+d-y)\bigr)\pmod{1024}.
\end{equation}
The factor after $256y$ satisfies
\[
 u\bigl(-1+qy+u(c+d-y)\bigr)\equiv
 \begin{cases}
 (q-1)u(u+1),&y\equiv1\pmod4,\\
 (q+1)u(u-1),&y\equiv-1\pmod4.
 \end{cases}
\]
Here $v_2(q)=0$, so again both expressions have valuation at least $2$.
This proves the required bound in the second case as well.  The convention
$v_2(0)=+\infty$ includes the axial points.
\end{proof}

\begin{lemma}
\label{lem:pointwise-two-adic}
For $z\in S^1\cap\Q(\ii)$, define
\begin{align*}
 \Phi(z)={}&-2\RePart z+2\ImPart z
 -7\RePart(z^3)+9\ImPart(z^3)\\
 &+2\RePart(z^4)+16\ImPart(z^4)
 +9\RePart(z^5)+7\ImPart(z^5).
\end{align*}
Then $\Phi(z)\in\Z_{(2)}$ and
\begin{equation}\label{eq:pointwise-two-adic}
 \Phi(z)\equiv2\pmod{2^{10}},\qquad v_2(\Phi(z)-2)\geq10.
\end{equation}
\end{lemma}

\begin{proof}
The coefficients come from a small choice of parameters in
Lemma~\ref{lem:four-param-lift}.  On the circle, the correction has the
Fourier expansion
\begin{equation}\label{eq:phi-parameter-fourier}
\begin{aligned}
 16\bigl(axy^2+bxy^4+cx^2y+dx^2y^3\bigr)
 ={}&(4a+2b)\RePart z+(4c+2d)\ImPart z\\
 &+(-4a-3b)\RePart(z^3)+(4c+d)\ImPart(z^3)\\
 &+b\RePart(z^5)-d\ImPart(z^5).
\end{aligned}
\end{equation}
This follows by expanding $\RePart(z^3)$, $\ImPart(z^3)$,
$\RePart(z^5)$, and $\ImPart(z^5)$ and reducing with $x^2+y^2=1$.

Set $h=q=1$ in the preceding proof.  Then $c=4$, $b=4-a$, and the
remaining conditions are $a\equiv-1\pmod4$ and $d\equiv1\pmod4$.
The real and imaginary coefficients become
\[
 (2a+8,-a-12,4-a),\qquad (16+2d,16+d,-d).
\]
Their squared Euclidean norms are, respectively,
\begin{align*}
 (2a+8)^2+(-a-12)^2+(4-a)^2&=6(a+4)^2+128,\\
 (16+2d)^2+(16+d)^2+d^2&=6(d+8)^2+128.
\end{align*}
The admissible choices closest to $-4$ and $-8$ are $a=-5$ and $d=-7$.
Taking $(a,b,c,d)=(-5,9,4,-7)$ therefore gives the triples $(-2,-7,9)$
and $(2,9,7)$.  Adding $\Phi_4$ recovers $\Phi$, and
Lemma~\ref{lem:four-param-lift} proves the congruence.  The norm
minimization here is restricted to $h=q=1$; we do not claim that the
choice is globally minimal or unique.
\end{proof}

\begin{proof}[Proof of Theorem~\ref{thm:fifth-divisibility}]
Let
\[
 P_k=\sum_{j=1}^N z_j^k.
\]
The assumed vanishing gives
\begin{align*}
 \sum_{j=1}^N\Phi(z_j)
 ={}&-2\RePart P_1+2\ImPart P_1
 -7\RePart P_3+9\ImPart P_3\\
 &+2\RePart P_4+16\ImPart P_4
 +9\RePart P_5+7\ImPart P_5=0.
\end{align*}
On the other hand, Lemma~\ref{lem:pointwise-two-adic} gives
\[
 0=\sum_{j=1}^N\Phi(z_j)\equiv2N\pmod{1024},
\]
so $v_2(N)=v_2(2N)-1\geq9$.  The fourth-moment theorem gives
$v_3(N)\geq2$, and hence $2^9 3^2=4608\mid N$.
\end{proof}

\subsection{Divisibility at strengths six and seven}
\label{subsec:sixth-divisibility}

\begin{theorem}\label{thm:sixth-divisibility}
Let $z_1,\ldots,z_N\in S^1\cap\Q(\ii)$, with repetitions allowed.  If
\[
 \sum_{j=1}^N z_j^k=0\qquad(1\leq k\leq6),
\]
then
\[
 9216\mid N.
\]
In particular, every rational equal-weight spherical $6$-design or
$7$-design on $S^1$ has cardinality divisible by $9216$.
\end{theorem}

A different linear combination of Fourier modes, denoted by $\Psi_6$,
improves the pointwise congruence by one more power of two.

\begin{lemma}
\label{lem:pointwise-six-two-adic}
For $z\in S^1\cap\Q(\ii)$, define
\begin{align*}
 \Psi_6(z)={}&4\ImPart z
 +\RePart(z^2)+\ImPart(z^2)
 -16\RePart(z^3)+2\ImPart(z^3)\\
 &+2\RePart(z^4)+16\ImPart(z^4)
 +16\RePart(z^5)-2\ImPart(z^5)\\
 &-\RePart(z^6)-11\ImPart(z^6).
\end{align*}
Then $\Psi_6(z)\in\Z_{(2)}$ and
\begin{equation}\label{eq:pointwise-six-two-adic}
 \Psi_6(z)\equiv2\pmod{2^{11}},\qquad v_2(\Psi_6(z)-2)\geq11.
\end{equation}
\end{lemma}

\begin{proof}
Write $z=x+\ii y$.  The two coordinate cases in
\eqref{eq:local-coordinate-valuations} are handled by different
factorizations.  Expanding and using $x^2+y^2=1$ gives
\begin{equation}\label{eq:psi-six-y-factor}
 \Psi_6(x+\ii y)-2=-32y^2(y-1)A_6(x,y),
\end{equation}
where
\[
 A_6(x,y)=11xy^2+3xy-4x-y^3+y,
\]
and also
\begin{equation}\label{eq:psi-six-x-factor}
 \Psi_6(x+\ii y)-2=-32x(x-1)(x+1)B_6(x,y),
\end{equation}
where
\[
 B_6(x,y)=x^3+11x^2y-8x^2+xy-x-4y+4.
\]

If $v_2(y)\geq2$, then $v_2(A_6(x,y))\geq2$ and $v_2(y-1)=0$.
The first factorization, \eqref{eq:psi-six-y-factor}, gives
\[
 v_2(\Psi_6(z)-2)
 \geq5+2v_2(y)+v_2(A_6(x,y))\geq11.
\]

In the remaining case, $v_2(x)\geq2$ and
$y,x-1,x+1\in\Z_{(2)}^\times$.  If $v_2(x)\geq3$, each of the first
five terms of $B_6(x,y)$ lies in $8\Z_{(2)}$, as does $4(1-y)$.
Thus $v_2(B_6(x,y))\geq3$, and
\[
 v_2(\Psi_6(z)-2)\geq5+v_2(x)+v_2(B_6(x,y))\geq11.
\]
If $v_2(x)=2$, write $x=4u$ with $u\in\Z_{(2)}^\times$.  Then
\[
 B_6(x,y)\equiv xy-x-4y+4=(x-4)(y-1)\pmod{16}.
\]
Now $v_2(x-4)=2+v_2(u-1)\geq3$ and $v_2(y-1)\geq1$, so
$v_2(B_6(x,y))\geq4$.  The second factorization,
\eqref{eq:psi-six-x-factor}, gives $v_2(\Psi_6(z)-2)\geq5+2+4=11$.
This proves \eqref{eq:pointwise-six-two-adic} in all cases.
\end{proof}

\begin{proof}[Proof of Theorem~\ref{thm:sixth-divisibility}]
Put $P_k=\sum_{j=1}^N z_j^k$.  The assumed vanishing gives
\begin{align*}
 \sum_{j=1}^N\Psi_6(z_j)={}&
 4\ImPart P_1+\RePart P_2+\ImPart P_2
 -16\RePart P_3+2\ImPart P_3\\
 &+2\RePart P_4+16\ImPart P_4
 +16\RePart P_5-2\ImPart P_5\\
 &-\RePart P_6-11\ImPart P_6=0.
\end{align*}
Lemma~\ref{lem:pointwise-six-two-adic} now yields
\[
 0=\sum_{j=1}^N\Psi_6(z_j)\equiv2N\pmod{2048},
\]
so $v_2(N)\geq10$.  Since $P_4=0$, the fourth-moment theorem also gives
$v_3(N)\geq2$.  Together these imply $2^{10}3^2=9216\mid N$.
\end{proof}

\subsection{Sharper divisibility under symmetry}\label{subsec:symmetric-divisibility}

Symmetry strengthens the divisibility obtained from the fourth moment
alone.  Start with a finite $K$-invariant rational set, written as
\begin{equation}\label{eq:K-decomp}
 Y=\varepsilon_1\{1,-1\}\sqcup
   \varepsilon_{\ii}\{\ii,-\ii\}\sqcup
   \bigsqcup_{r=1}^M\Korb(w_r),
 \qquad \varepsilon_1,\varepsilon_{\ii}\in\{0,1\},
\end{equation}
where the four-point orbits are ordinary.

\begin{theorem}
\label{thm:K-intro}\label{thm:K-divisibility}
If a finite rational spherical $4$-design $Y\subset S^1$ is $K$-invariant,
then
\[
 4608\mid |Y|.
\]
Every such design is automatically a spherical $5$-design.  Thus, within
the $K$-invariant class, the sharp divisor $4608$ follows from the fourth
moment alone; in particular, $X_5$ also attains this symmetry-class bound.
\end{theorem}

\begin{proof}
For the representatives in \eqref{eq:K-decomp}, choose a common
denominator $C$ by Lemma~\ref{lem:common-denominator} and write
\[
 w_r=\frac{A_r+\ii B_r}{C},
 \qquad q_r=\frac{A_rB_r}{12}\in\Z.
\]
Each exceptional two-point orbit contributes $2$ to $P_4$.
By Lemma~\ref{lem:orbit-moments}, an ordinary orbit contributes
\[
 4\left(1-\frac{1152q_r^2}{C^4}\right).
\]
The equation $P_4(Y)=0$ gives
\[
 0=2(\varepsilon_1+\varepsilon_{\ii})+4M
   -\frac{4608}{C^4}\sum_{r=1}^Mq_r^2.
\]
The first two terms sum to $|Y|$, so
\[
 |Y|C^4=4608\sum_{r=1}^Mq_r^2.
\]
Because $v_2(C)=v_3(C)=0$, we have $v_2(|Y|)\geq9$ and
$v_3(|Y|)\geq2$, which proves $4608\mid |Y|$.
Corollary~\ref{cor:automatic} gives the $5$-design property, and the
example in Theorem~\ref{thm:main} attains the bound.
\end{proof}

\begin{theorem}
\label{thm:D4-divisibility}
If $Y\subset S^1\cap\Q^2$ is a finite $D_4$-invariant spherical $4$-design,
then
\[
 9216\mid |Y|.
\]
Moreover, the axis orbit $\muFour$ does not occur in $Y$.
\end{theorem}

\begin{proof}
Every such rational set decomposes as
\[
 Y=\varepsilon\muFour\sqcup\bigsqcup_{r=1}^M\Orb(w_r),
 \qquad \varepsilon\in\{0,1\}.
\]
The diagonal points have coordinates $\pm1/\sqrt2$ and are not rational.
Choose a common denominator $C$ using Lemma~\ref{lem:common-denominator},
and write
\[
 w_r=\frac{A_r+\ii B_r}{C},
 \qquad q_r=\frac{A_rB_r}{12}\in\Z.
\]
The fourth moment gives
\[
 (2M+\varepsilon)C^4=2304\sum_{r=1}^Mq_r^2.
\]
Since $C\in\Z_{(2)}^\times\cap\Z_{(3)}^\times$, it follows that
$v_2(2M+\varepsilon)\geq8$ and $v_3(2M+\varepsilon)\geq2$.
Reduction modulo $2$ forces $\varepsilon=0$.  We are left with
$v_2(M)\geq7$ and $v_3(M)\geq2$, hence $9216\mid |Y|=8M$.
\end{proof}

\begin{theorem}\label{thm:eighth-obstruction}
Every $D_4$-invariant rational spherical $8$-design has cardinality divisible
by
\[
 7\,077\,888.
\]
\end{theorem}

\begin{proof}
By Theorem~\ref{thm:D4-divisibility}, the axis orbit is absent.
The design is therefore a union of $M$ ordinary orbits.  Choose $C$ as
in Lemma~\ref{lem:common-denominator} and write
\[
 w_r=\frac{A_r+\ii B_r}{C},
 \qquad q_r=\frac{A_rB_r}{12}\in\Z,
\]
and set
\[
 x_r=\RePart(w_r^4)=1-\frac{1152q_r^2}{C^4}.
\]
The fourth and eighth moments give
\[
 \sum_rx_r=0,
 \qquad
 \sum_r(2x_r^2-1)=0.
\]
The first identity gives
\[
 \sum_rq_r^2=\frac{MC^4}{1152}.
\]
Substituting this into the second yields
\[
 2\cdot1152^2\sum_rq_r^4=3MC^8,
\]
or equivalently
\[
 2^{15}3^3\sum_rq_r^4=MC^8.
\]
Taking valuations and using $v_2(C)=v_3(C)=0$ gives $v_2(M)\geq15$
and $v_3(M)\geq3$.  Since $|Y|=8M$, we obtain
$2^{18}3^3=7\,077\,888\mid |Y|$.
\end{proof}

\subsection{Universal local divisibility in arbitrary strength}
\label{subsec:universal-Dt-proof}

The general divisor combines the bounds already proved with congruences
at the inert primes.  For $t\geq1$, put
\[
 r=\left\lfloor\frac t4\right\rfloor,
 \qquad B_s=\binom{2s}{s}\quad(0\leq s\leq r),
\]
and define the low-degree $2$-adic seed
\begin{equation}\label{eq:s2-seed}
 s_2(t)=
 \begin{cases}
  0,&t\leq3,\\
  8,&t=4,\\
  9,&t=5,\\
  10,&t\geq6.
 \end{cases}
\end{equation}
Set
\begin{equation}\label{eq:alpha23-definition}
 \alpha_2(t)=\max\left\{s_2(t),\,
   \max_{0\leq s\leq r}\bigl(8s-v_2(B_s)\bigr)\right\},
 \qquad
 \alpha_3(t)=\max_{0\leq s\leq r}\bigl(2s-v_3(B_s)\bigr)
\end{equation}
and define
\begin{equation}\label{eq:Dt-definition}
 D_t=
 2^{\alpha_2(t)}3^{\alpha_3(t)}
 \prod_{\substack{p\in\mathcal P_-\\7\leq p\leq t-1}}
 p^{\lfloor t/(p+1)\rfloor}.
\end{equation}
The product is finite.  The central binomial valuations are given by
Legendre's formula:
\[
 v_p(B_s)=\sum_{\ell\geq1}
 \left(\left\lfloor\frac{2s}{p^\ell}\right\rfloor
       -2\left\lfloor\frac{s}{p^\ell}\right\rfloor\right).
\]
Indeed, $v_p(n!)=\sum_{\ell\geq1}\lfloor n/p^\ell\rfloor$ counts the
multiples of each prime power among $1,\ldots,n$; the sum is finite.

\begin{theorem}
\label{thm:universal-Dt}
Let $z_1,\ldots,z_N\in S^1\cap\Q(\ii)$, with repetitions allowed, and suppose
\[
 \sum_{j=1}^N z_j^k=0\qquad(1\leq k\leq t).
\]
Then
\[
 D_t\mid N.
\]
Consequently, whenever $\NQ(t)<\infty$, one has $D_t\mid\NQ(t)$.
\end{theorem}

\begin{corollary}
\label{cor:Dt-values-asymptotic}
The universal divisor recovers all sharp cardinalities proved above:
\[
 D_4=2304,
 \qquad D_5=4608,
 \qquad D_6=D_7=9216.
\]
Moreover,
\[
 D_t=2^{15}3^3\cdot7=6\,193\,152
 \qquad(8\leq t\leq11),
\]
and
\[
 D_t=2^{22}3^6\cdot7\cdot11=235\,438\,866\,432
 \qquad(12\leq t\leq15).
\]
As $t\to\infty$,
\begin{equation}\label{eq:Dt-asymptotic}
 \log D_t=\left(\frac12+o(1)\right)t\log t.
\end{equation}
Thus, at every strength for which a finite rational equal-weight design
exists,
\[
 \NQ(t)\geq t^{(1/2+o(1))t}.
\]
\end{corollary}

\begin{remark}
\label{rem:Dt-not-maximal}
Including the seed $s_2(t)$ preserves the sharp low-degree congruences
and gives $D_t=\NQ(t)$ for $t=4,5,6,7$.  For $t\geq8$, we claim only
that $D_t$ is a necessary divisor, not that it is the largest one.
Additional congruences at $p=2$ or $p\in\mathcal P_-$ can be sought on
the finite set
\[
 U_{p,m}=\{(x,y)\in(\Z/p^m\Z)^2:x^2+y^2=1\}.
\]
Every rational circle point reduces to $U_{p,m}$ by
\eqref{eq:local-coordinate-valuations}.  Evaluating $1,\Re z^k,\Im z^k$
on this set gives a finite matrix.  Smith normal form reduces the resulting
linear congruence problem; see Cohen~\cite[Chapter~2]{Cohen93} for the
standard integer-matrix algorithms.  This calculation concerns necessary
congruences and does not assert that a residue solution lifts to a
rational design.
\end{remark}

For Theorem~\ref{thm:universal-Dt}, the primes $2$ and $3$ are treated
using powers of $\RePart(z^4)-1$.  At the other inert primes, we use the
norm-one relation and Frobenius.

\begin{lemma}
\label{lem:fourth-character-binomial}
Under the hypotheses of Theorem~\ref{thm:universal-Dt}, for every
integer $s$ with $0\leq s\leq\lfloor t/4\rfloor$ and
$B_s=\binom{2s}{s}$, one has
\[
 v_2(N)\geq8s-v_2(B_s),
 \qquad
 v_3(N)\geq2s-v_3(B_s).
\]
\end{lemma}

\begin{proof}
For $z=x+\ii y$, \eqref{eq:local-coordinate-valuations} gives
\[
 v_2(x^2y^2)\geq4,\qquad v_3(x^2y^2)\geq2.
\]

On $S^1$,
\begin{equation}\label{eq:fourth-character-factor}
 \RePart(z^4)-1=-8x^2y^2.
\end{equation}
Therefore the degree-$4s$ polynomial
\begin{equation}\label{eq:Fs-definition}
 F_s(z)=\bigl(\RePart(z^4)-1\bigr)^s
\end{equation}
satisfies, pointwise,
\[
 v_2(F_s(z))\geq7s,\qquad v_3(F_s(z))\geq2s.
\]
Since its degree is at most $t$, its design average equals its circular
average.  The identity $\cos4\theta-1=-2\sin^2(2\theta)$ gives
\begin{equation}\label{eq:Fs-average}
 \frac1{2\pi}\int_0^{2\pi}(\cos4\theta-1)^s\,d\theta
 =(-1)^s\frac{B_s}{2^s}.
\end{equation}
Consequently
\[
 \sum_{j=1}^N F_s(z_j)=(-1)^sN\frac{B_s}{2^s}.
\]
Taking valuations and using
$v_p(\sum_jF_s(z_j))\geq\min_j v_p(F_s(z_j))$, we obtain
\[
 v_2(N)+v_2(B_s)-s\geq7s,
 \qquad
 v_3(N)+v_3(B_s)\geq2s,
\]
which is the claim.
\end{proof}

\begin{lemma}
\label{lem:inert-prime-divisibility}
Let $p\in\mathcal P_-$ with $p\geq7$.  Under the hypotheses of
Theorem~\ref{thm:universal-Dt},
\[
 p^{\lfloor t/(p+1)\rfloor}\mid N.
\]
\end{lemma}

\begin{proof}
Write $r_p=\lfloor t/(p+1)\rfloor$; only $r_p\geq1$ needs proof.
By \eqref{eq:local-coordinate-valuations}, each $z_j$ belongs to
$\mathcal O_{(p)}=\Z_{(p)}[\ii]$ and has norm $1$.
Since $\bigl(\frac{-1}{p}\bigr)=-1$, the polynomial $T^2+1$ is
irreducible over $\mathbb F_p$, so
\[
 \mathcal O_{(p)}/p\mathcal O_{(p)}
 \cong\mathbb F_p[T]/(T^2+1)\cong\mathbb F_{p^2}.
\]
On this residue field, Frobenius is conjugation.  Thus the reduction
$\zeta_j$ of $z_j$ satisfies
\[
 \Norm_{\mathbb F_{p^2}/\mathbb F_p}(\zeta_j)
 =\zeta_j^{p+1}=1.
\]
Thus $(z_j^{p+1}-1)^{r_p}\in p^{r_p}\mathcal O_{(p)}$.  Expanding and
summing, with $P_0=N$, gives
\begin{align*}
 \sum_{j=1}^N(z_j^{p+1}-1)^{r_p}
 &=\sum_{\ell=0}^{r_p}(-1)^{r_p-\ell}
   \binom{r_p}{\ell}P_{\ell(p+1)}\\
 &=(-1)^{r_p}N.
\end{align*}
The terms with $\ell\geq1$ vanish because $\ell(p+1)\leq t$.  Hence
\[
 N\in p^{r_p}\mathcal O_{(p)}\cap\Z=p^{r_p}\Z,
\]
which proves the required divisibility.
\end{proof}

\begin{proof}[Proof of Theorem~\ref{thm:universal-Dt}]
Apply Lemma~\ref{lem:fourth-character-binomial} over $0\leq s\leq r$.
Theorems~\ref{thm:universal-divisibility}, \ref{thm:fifth-divisibility},
and~\ref{thm:sixth-divisibility} supply the seed $s_2(t)$: they give
$v_2(N)$ at least $8$, $9$, and $10$ at strengths at least $4$, $5$,
and $6$, respectively.  Thus
\[
 v_2(N)\geq\alpha_2(t),\qquad v_3(N)\geq\alpha_3(t).
\]
For each $p\in\mathcal P_-$ with $p\geq7$,
Lemma~\ref{lem:inert-prime-divisibility} gives
$v_p(N)\geq\lfloor t/(p+1)\rfloor$.  We have therefore shown
$v_p(N)\geq v_p(D_t)$ at every prime, which proves $D_t\mid N$.
\end{proof}

\begin{proof}[Proof of Corollary~\ref{cor:Dt-values-asymptotic}]
When $4\leq t\leq7$, we have $r=1$ and $B_r=2$.
Using the seeds in \eqref{eq:s2-seed} gives $D_4=2304$, $D_5=4608$,
and $D_6=D_7=9216$.  For $8\leq t\leq11$, we have $r=2$ and $B_r=6$.
The only additional inert prime is $7$, so
\[
 D_t=2^{15}3^3\cdot7=6\,193\,152.
\]
For $12\leq t\leq15$, the values are $r=3$ and $B_r=20$.
The additional inert primes are $7$ and $11$, giving
\[
 D_t=2^{22}3^6\cdot7\cdot11=235\,438\,866\,432.
\]

The powers of $2$ and $3$ contribute $O(t)$ to $\log D_t$.
For the remaining primes, apply the prime number theorem in the
progression $3\pmod4$, or equivalently for $\chi_{-4}(p)=-1$
\cite{Davenport}.  Partial summation gives
\[
 \sum_{\substack{p\in\mathcal P_-\\p\leq x}}\frac{\log p}{p}
 =\frac12\log x+O(1).
\]
Hence
\begin{align*}
 \sum_{\substack{p\in\mathcal P_-\\7\leq p\leq t-1}}
 \left\lfloor\frac{t}{p+1}\right\rfloor\log p
 &=t\sum_{\substack{p\in\mathcal P_-\\p\leq t-1}}
   \frac{\log p}{p+1}+O(t)\\
 &=\frac12t\log t+O(t).
\end{align*}
This proves \eqref{eq:Dt-asymptotic}.  Whenever $\NQ(t)$ is finite,
Theorem~\ref{thm:universal-Dt} then gives the stated lower bound.
\end{proof}

\section{Constructions of optimal designs at strengths 5, 6, and 7}
\label{sec:optimal567}\label{sec:designs}

We construct $(X_5,X_7)$ and $(\widetilde X_5,\widetilde X_7)$ on a
single Gaussian shell.  The construction has two steps: select orbits
whose fourth moments sum to zero, then orient their four-point suborbits
so that the second moments also sum to zero.  An exact orbit trade and a
new orientation give the second pair.  The $4608$-point sets have exact
strength $5$; their $9216$-point quarter-turn completions have exact
strength $7$.  The divisibility results of
Section~\ref{sec:fourier-divisibility} prove optimality at strengths $5$,
$6$, and $7$ in the usual nested sense.  We distinguish the two pairs
modulo $O(2)$ by their absolute sixth and eighth moments.

\subsection{The common Gaussian shell}\label{sec:H5-shell}

Both $X_5$ and its quarter-turn completion $X_7$ use the shell with
hypotenuse $\Hfive$.  We index candidate representatives by $i$ and
selected representatives by $j$.  Put
\[
 (\mathfrak p_1,\ldots,\mathfrak p_7)
 =(1+2\ii,2+3\ii,1+4\ii,2+5\ii,1+6\ii,4+5\ii,2+7\ii),
\]
\[
 m_1=2,\qquad m_2=\cdots=m_7=1.
\]
Their norms are the distinct rational primes $5,13,17,29,37,41,53$.
Each $\mathfrak p_j$ is therefore a Gaussian prime.  Since its
coordinates are nonzero and unequal in absolute value,
$\mathfrak p_j$ is not associated to $\overline{\mathfrak p_j}$.
Thus the fourteen primes $\mathfrak p_j,\overline{\mathfrak p_j}$ are
pairwise nonassociate.  Define
\begin{equation}\label{eq:H5}
 \Hfive=\prod_{j=1}^7\Norm(\mathfrak p_j)^{m_j}
 =5^2\cdot13\cdot17\cdot29\cdot37\cdot41\cdot53
 =12\,882\,250\,225.
\end{equation}

\begin{remark}\label{rem:shell-choice}
If $\mathfrak p=u+\ii v$ has rational prime norm $p\equiv1\pmod4$, then
\[
 \frac{\mathfrak p}{\overline{\mathfrak p}}
 =\frac{u^2-v^2}{p}+\ii\frac{2uv}{p}
 \in S^1\cap\Q(\ii).
\]
This gives a rational rotation for each split prime.  We use the seven
smallest primes congruent to $1$ modulo $4$.

For positive base multiplicities $m_1,\ldots,m_r$, the exponent box
$-m_j\leq e_j\leq m_j$ has
\[
 \frac{\prod_j(2m_j+1)-1}{2}
\]
nonzero classes under $e\sim-e$.  Proposition~\ref{prop:candidates}
identifies them with ordinary $D_4$-orbits for our choice of Gaussian
primes.  Taking all seven multiplicities equal to $1$ gives only
\[
 \frac{3^7-1}{2}=1093<1152
\]
candidates.  Increasing one multiplicity to $2$ gives
\[
 \frac{5\cdot3^6-1}{2}=1822>1152.
\]
Increasing the multiplicity at norm $5$, rather than at one of the other
six norms, gives the smallest common denominator among these choices.
We therefore take
\[
 (m_1,\ldots,m_7)=(2,1,1,1,1,1,1).
\]
The box now contains enough candidates to select $1152$ ordinary orbits.
It remains to find a selection and orientation satisfying the moment
equations.
\end{remark}

Define
\[
 \Efive=\left\{e=(e_1,\ldots,e_7)\in\Z^7:
 \begin{array}{l}
 -2\leq e_1\leq2,\quad -1\leq e_j\leq1\ (2\leq j\leq7),\\
 e\neq0,\quad\text{the first nonzero coordinate of $e$ is positive}
 \end{array}\right\}.
\]
For $e\in\Efive$, put
\begin{equation}\label{eq:gamma5}
 \gamma_e=
 \prod_{j=1}^7\mathfrak p_j^{m_j+e_j}\overline{\mathfrak p_j}^{\,m_j-e_j}
 =A_e+\ii B_e.
\end{equation}
Then
\begin{equation}\label{eq:ratio-and-shell}
 \frac{\gamma_e}{\Hfive}
 =\prod_{j=1}^7\left(\frac{\mathfrak p_j}{\overline{\mathfrak p_j}}\right)^{e_j}
 \in S^1\cap\Q(\ii),
 \qquad
 A_e^2+B_e^2=\Hfive^2,
 \qquad \Delta_e=A_e^2-B_e^2.
\end{equation}

\begin{proposition}\label{prop:candidates}
The set $\Efive$ has $1822$ elements.  Every $\Orb(\gamma_e/\Hfive)$ has eight
points, and these orbits are pairwise disjoint.
\end{proposition}

\begin{proof}
The box contains $5\cdot3^6=3645$ vectors.  Remove zero and choose one
vector from each pair $\{e,-e\}$; this leaves $1822$ representatives.

An orbit can have fewer than eight points only on an axis or a diagonal.
Suppose $A_eB_e=0$.  Then $\gamma_e=u\Hfive$ for some unit
$u\in\{\pm1,\pm\ii\}$.  Comparing exponents of the fourteen pairwise
nonassociate Gaussian primes in
\[
 \gamma_e=
 \prod_{j=1}^7\mathfrak p_j^{m_j+e_j}
 \overline{\mathfrak p_j}^{\,m_j-e_j},
 \qquad
 \Hfive=\prod_{j=1}^7
 \mathfrak p_j^{m_j}\overline{\mathfrak p_j}^{\,m_j},
\]
forces $m_j+e_j=m_j$ for every $j$, hence $e=0$, a contradiction.
The diagonal case $|A_e|=|B_e|$ would imply $\Hfive^2=2A_e^2$, also
impossible because $\Hfive$ is odd.  All candidate orbits have eight points.

If two candidate orbits meet, then for some Gaussian unit $u$ their
representatives satisfy

\[
 \gamma_f=u\gamma_e
 \qquad\text{or}\qquad
 \gamma_f=u\overline{\gamma_e}.
\]
In the first case, comparison of the $\mathfrak p_j$-adic valuations gives
\[
 m_j+f_j=m_j+e_j\qquad(1\leq j\leq7),
\]
and hence $f=e$.  In the second case,
\[
 m_j+f_j=m_j-e_j\qquad(1\leq j\leq7),
\]
and hence $f=-e$.  This is excluded by the convention that the first
nonzero coordinate is positive.  Thus two candidate orbits meet only
when $f=e$.
\end{proof}

Put
\[
 q_e=\frac{A_eB_e}{12}\in\Z.
\]

\subsection{Selection and orientation}\label{sec:selection}

The first step is to select $1152$ of the $1822$ ordinary $D_4$-orbits
with total fourth moment zero.  Index the candidates from
Subsection~\ref{sec:H5-shell} by $\mathcal I_5=\{0,\ldots,1821\}$ and
write $q_i=A_iB_i/12$.  Equation~\eqref{eq:Rew4} reduces the selection
to a fixed-cardinality subset sum:
\begin{equation}\label{eq:search-selection-model}
 b_i\in\{0,1\},\qquad
 \sum_{i\in\mathcal I_5}b_i=1152,\qquad
 \sum_{i\in\mathcal I_5}b_iq_i^2=\Hfive^4.
\end{equation}
Next choose one of the two $K$-orbits in each selected orbit.
Multiplication by $\ii$ negates the square of a representative but leaves
its fourth power unchanged.  Thus the orientation problem on the selected
support $\mathcal S$ is
\[
 \sigma_i\in\{-1,1\},\qquad
 \sum_{i\in\mathcal S}\sigma_i|\Delta_i|=0.
\]
Algorithm~\ref{alg:five} searches for solutions by capped,
problem-specific variants of Wagner's generalized-birthday list-merging
method~\cite{Wagner02}.
The candidates are shuffled into $128$ blocks of size $14$ or $15$,
with nine indices selected from each block.  For the orientation
half-sum, the selected entries are arranged in $128$ blocks of nine.
Every subset of an orientation block is allowed, so the number of
positive signs is not fixed.

The block restrictions and list caps make the search incomplete.
If orientation fails, the run stops; it does not return to an earlier
selection.  A new selection therefore requires a new invocation.
The proof uses the fixed certificate below, independently of the search.
The implemented runs found different designs, recorded in
Appendix~\ref{app:search-provenance}.  No additional search is needed
for $X_7$: it is the quarter-turn completion of $X_5$.

\begin{proposition}\label{prop:selection-certificate}
There are distinct exponent vectors
\[
 e_0,\ldots,e_{1151}\in\Efive
\]
and signs $\sigma_j\in\{\pm1\}$ with the following properties.  Write
\[
 \gamma_j=\gamma_{e_j}=A_j+\ii B_j,
 \qquad \Delta_j=A_j^2-B_j^2,
 \qquad q_j=\frac{A_jB_j}{12}\in\Z.
\]
Then
\begin{equation}\label{eq:first-certificate}
 \sum_{j=0}^{1151}q_j^2=\Hfive^4,
\end{equation}
and
\begin{equation}\label{eq:second-certificate}
 \sum_{j=0}^{1151}\sigma_j|\Delta_j|=0.
\end{equation}
For each $j$, choose
$\beta_j=\widehat A_j+\ii\widehat B_j\in\{\gamma_j,\ii\gamma_j\}$
so that
\begin{equation}\label{eq:beta2}
 \RePart(\beta_j^2)=\widehat A_j^{\,2}-\widehat B_j^{\,2}
 =\sigma_j|\Delta_j|.
\end{equation}
The resulting collection also satisfies
\begin{equation}\label{eq:halving-parity}
 \sum_{j=0}^{1151}\frac{\widehat A_j\widehat B_j}{12}\equiv1\pmod2.
\end{equation}
Finally, define
\begin{equation}\label{eq:X5}
 X_5=\bigsqcup_{j=0}^{1151}\Korb\left(\frac{\beta_j}{\Hfive}\right).
\end{equation}
Its sixth moment is the exact nonzero integer certificate
{\small
\begin{equation}\label{eq:sixth-certificate}
 \Hfive^6P_6(X_5)
 =-462\,563\,147\,140\,280\,729\,960\,096\,119\,277\,432\,099\,780\,213\,746\,861\,044\,063\,106\,513\,920
 \neq0.
\end{equation}}
\end{proposition}

\begin{proof}
The selection and orientation data are listed in
Appendix~\ref{app:certificate-pair}.  Reconstruct the Gaussian integers
from these data and sum over their four-point orbits.  The displayed
identities, including \eqref{eq:sixth-certificate}, follow by exact
integer arithmetic.
\end{proof}

\subsection{The first optimal pair}\label{subsec:first-optimal-pair}

Proposition~\ref{prop:candidates} makes the selected $D_4$-orbits
pairwise disjoint.  Their four-point $K$-suborbits defining $X_5$ in
\eqref{eq:X5} are consequently disjoint as well.

\begin{proof}[Proof of Theorem~\ref{thm:main}]
The set has
\[
 |X_5|=4\cdot1152=4608,
\]
and lies in $\Hfive^{-1}\Z^2\cap S^1$.  Antipodality gives vanishing
of the odd moments.  For the second moment, combine \eqref{eq:beta2}
with \eqref{eq:second-certificate}:
\[
 P_2(X_5)=\frac4{\Hfive^2}\sum_{j=0}^{1151}\RePart(\beta_j^2)=0.
\]
For the fourth moment, use $\beta_j^4=\gamma_j^4$ and
$A_j^2+B_j^2=\Hfive^2$.  We have
\[
 \RePart(\gamma_j^4)
 =(A_j^2+B_j^2)^2-8A_j^2B_j^2
 =\Hfive^4-1152q_j^2.
\]
Equation~\eqref{eq:first-certificate} then yields
\[
 P_4(X_5)
 =\frac4{\Hfive^4}\sum_{j=0}^{1151}\RePart(\beta_j^4)
 =0.
\]
We have $P_k(X_5)=0$ for $1\leq k\leq5$, so
Lemma~\ref{lem:fourier} proves that $X_5$ is a $5$-design.
The nonzero sixth moment in \eqref{eq:sixth-certificate} makes the
strength exactly $5$.  Its cardinality is optimal by
Theorem~\ref{thm:fifth-divisibility}.
\end{proof}

\begin{proof}[Proof of Corollary~\ref{cor:NQ5}]
The lower bound $\NQ(5)\geq4608$ is
Theorem~\ref{thm:fifth-divisibility}; the $4608$-point design of
Theorem~\ref{thm:main} attains it.
\end{proof}

Define
\begin{equation}\label{eq:X7}
 X_7=X_5\sqcup\ii X_5
 =\bigsqcup_{j=0}^{1151}\Orb\left(\frac{\gamma_j}{\Hfive}\right).
\end{equation}
Within each selected $D_4$-orbit, the $K$-orbit in $X_5$ and its
quarter-turn are disjoint and together fill the orbit.  Thus $|X_7|=9216$.

\begin{proof}[Proof of Theorem~\ref{thm:parent-intro}]
The $D_4$ symmetry of $X_7$ leaves only the fourth moment to check
among degrees one through $7$.  It vanishes because
\[
 P_4(X_7)=P_4(X_5)+P_4(\ii X_5)=2P_4(X_5)=0.
\]
Hence $X_7$ is a spherical $7$-design.  Were it an $8$-design,
Theorem~\ref{thm:eighth-obstruction} would force
\[
 7\,077\,888\mid9216,
\]
which is impossible.  The strength is therefore exactly $7$.
Optimality at both strengths follows from
Theorem~\ref{thm:sixth-divisibility}.  The assertion for the
$D_4$-invariant class follows from Theorem~\ref{thm:D4-divisibility}.
\end{proof}

\begin{proof}[Proof of Corollary~\ref{cor:NQ67}]
The lower bound at both strengths is $9216$ by
Theorem~\ref{thm:sixth-divisibility}.  The $7$-design in
Theorem~\ref{thm:parent-intro} has this cardinality and is also a
$6$-design.
\end{proof}

\subsection{A second optimal parent from an exact orbit trade}
\label{subsec:second-x7}

Retain the zero-based lexicographic indices
$\mathcal I_5=\{0,\ldots,1821\}$ for the $1822$ candidate $D_4$-orbits.
For the trade, replace each representative by the one with positive
coordinates in increasing order, without changing its index.
For $i\in\mathcal I_5$, write
\[
 \gamma_i=A_i+\ii B_i,
 \qquad 0<A_i<B_i,
 \qquad \Delta_i=A_i^2-B_i^2,
 \qquad q_i=\frac{A_iB_i}{12}\in\Z,
\]
and put $\mathcal O_i=\Orb(\gamma_i/\Hfive)$.  Let
$\mathcal S\subset\mathcal I_5$ be the $1152$-element selection defining
$X_7$.

A selection $\mathcal S'$ with the same cardinality and $q_i^2$-sum
as $\mathcal S$ gives a trade by removing
$\mathcal R=\mathcal S\setminus\mathcal S'$ and adding
$\mathcal A=\mathcal S'\setminus\mathcal S$.
Algorithm~\ref{alg:second-seven} searches for $\mathcal S'$ by running
the $128$-block search with a different seed.  It then checks the trade
identities and the absolute eighth moment; it does not prescribe the
size of the trade.  The certificate used here exchanges $326$ orbits
in each direction.  The separately recorded appendix pair exchanges
$435$.

\begin{proposition}\label{prop:x7-trade326}
There are sets $\mathcal R\subset\mathcal S$ and
$\mathcal A\subset\mathcal I_5\setminus\mathcal S$ such that
$|\mathcal R|=|\mathcal A|=326$ and
\begin{equation}\label{eq:x7-trade326}
\begin{aligned}
 \sum_{i\in\mathcal R}q_i^2
 &=\sum_{i\in\mathcal A}q_i^2\\
 &=7\,785\,922\,792\,245\,261\,434\,720\,277\,359\,475\,660\,903\,110.
\end{aligned}
\end{equation}
Put
\[
 \widetilde{\mathcal S}=(\mathcal S\setminus\mathcal R)\cup\mathcal A,
 \qquad
 \widetilde{X}_7=\bigsqcup_{i\in\widetilde{\mathcal S}}\mathcal O_i.
\]
Then $\widetilde{X}_7\subset \Hfive^{-1}\Z^2\cap S^1$ consists of $9216$
distinct points, is a spherical $7$-design of exact strength $7$, and
\[
 |X_7\mathbin{\triangle}\widetilde{X}_7|=5216.
\]
Moreover, its exact eighth-moment certificate satisfies
\begin{equation}\label{eq:x7tilde-p8-mod13}
 \Hfive^8P_8(\widetilde{X}_7)\equiv7\pmod{13},
\end{equation}
so in particular $P_8(\widetilde{X}_7)\neq0$.
\end{proposition}

\begin{proof}
Use the index sets $\mathcal R$ and $\mathcal A$ listed in
Appendix~\ref{app:certificate-traded-parent}.  Direct summation of the
shell data gives \eqref{eq:x7-trade326}.  Checking
$\mathcal R\subset\mathcal S$, $\mathcal A\cap\mathcal S=\varnothing$,
and $\widetilde{\mathcal S}=(\mathcal S\setminus\mathcal R)\cup\mathcal A$
then gives $|\widetilde{\mathcal S}|=1152$.
The selected orbits are disjoint by Proposition~\ref{prop:candidates},
so $\widetilde{X}_7$ has $8\cdot1152=9216$ distinct rational points.
The removed and added orbits account for the symmetric difference
$|X_7\mathbin{\triangle}\widetilde{X}_7|=8(|\mathcal R|+|\mathcal A|)=5216$.

For one candidate orbit,
\[
 P_4(\mathcal O_i)
 =\frac8{\Hfive^4}\RePart(\gamma_i^4)
 =\frac8{\Hfive^4}\bigl(\Hfive^4-1152q_i^2\bigr).
\]
For the original selection, \eqref{eq:first-certificate} gives
$\sum_{i\in\mathcal S}q_i^2=\Hfive^4$.  Since the trade preserves the
sum, we also have $\sum_{i\in\widetilde{\mathcal S}}q_i^2=\Hfive^4$.
It follows that
\[
 P_4(\widetilde{X}_7)
 =\frac8{\Hfive^4}\left(1152\Hfive^4
 -1152\sum_{i\in\widetilde{\mathcal S}}q_i^2\right)=0.
\]
Quarter-turn invariance gives vanishing in all degrees not divisible by
$4$.  Thus $P_k(\widetilde{X}_7)=0$ for $1\leq k\leq7$, and
Lemma~\ref{lem:fourier} proves the $7$-design property.
Reconstructing the same points and evaluating
$\Hfive^8P_8(\widetilde{X}_7)$ gives \eqref{eq:x7tilde-p8-mod13}.
The eighth moment is nonzero, so the strength is exactly $7$.
\end{proof}

\subsection{An optimal five-design in the traded parent}
\label{subsec:second-x5}

We must orient the traded support anew, indexing its signs by their
ranks in the new selection.  Algorithm~\ref{alg:second-five} uses the
same $128$-block half-sum search.  The $576$--$576$ sign balance of the
certificate below is not a search constraint: the two orientations
recorded in the appendix have $587$ and $574$ positive signs.

\begin{proposition}
\label{prop:x5tilde-orientation}
List the traded selection increasingly as
$\widetilde{\mathcal S}=\{i_0<\cdots<i_{1151}\}$.  There are signs
$\widetilde\sigma_r\in\{\pm1\}$, with exactly $576$ signs of each type, such
that
\begin{equation}\label{eq:x5tilde-orientation}
 \sum_{r=0}^{1151}\widetilde\sigma_r
 |\Delta_{i_r}|=0.
\end{equation}
Choose
\[
 \widetilde\beta_r\in
 \{\gamma_{i_r},\ \ii\gamma_{i_r}\}
\]
so that
\[
 \RePart(\widetilde\beta_r^2)
 =\widetilde\sigma_r|\Delta_{i_r}|,
\]
and define
\begin{equation}\label{eq:X5tilde}
 \widetilde X_5
 =\bigsqcup_{r=0}^{1151}
   \Korb\left(\frac{\widetilde\beta_r}{\Hfive}\right).
\end{equation}
Then $\widetilde X_5$ consists of $4608$ distinct rational points,
\[
 \widetilde X_7=\widetilde X_5\sqcup\ii\widetilde X_5,
\]
and $\widetilde X_5$ is a spherical $5$-design of exact strength $5$.
Its cleared sixth moment is
{\small
\begin{equation}\label{eq:x5tilde-p6}
 \Hfive^6P_6(\widetilde X_5)
 =1\,133\,586\,542\,449\,770\,272\,769\,317\,889\,152\,349\,343\,934\,086\,361\,589\,179\,131\,806\,720
 \neq0.
\end{equation}}
In particular, $\widetilde X_5$ is globally cardinality-optimal.
\end{proposition}

\begin{proof}
Apply the orientation data in Appendix~\ref{app:certificate-traded-child}
to the selection $\widetilde{\mathcal S}$ from
Proposition~\ref{prop:x7-trade326}.  Reconstruction verifies the sign
balance, \eqref{eq:x5tilde-orientation}, and distinctness of the points.
Their quarter-turn completion is $\widetilde X_7$.

The odd moments vanish by antipodality, and
\eqref{eq:x5tilde-orientation} gives
\[
 P_2(\widetilde X_5)
 =\frac4{\Hfive^2}\sum_{r=0}^{1151}
   \RePart(\widetilde\beta_r^2)=0.
\]
Multiplication by $\ii$ does not change a fourth power, so
\[
 P_4(\widetilde X_5)
 =\frac4{\Hfive^4}\sum_{r=0}^{1151}
   \RePart(\gamma_{i_r}^4)
 =\frac12P_4(\widetilde X_7)=0.
\]
Hence all moments through degree $5$ vanish.  In each selected
$D_4$-orbit, $\Korb(\widetilde\beta_r/\Hfive)$ and
$\Korb(\ii\widetilde\beta_r/\Hfive)$ are disjoint and exhaust the orbit,
so $\widetilde X_7=\widetilde X_5\sqcup\ii\widetilde X_5$.
The sixth moment computed from the same data is \eqref{eq:x5tilde-p6},
which proves exact strength $5$.  Corollary~\ref{cor:NQ5} gives
cardinality optimality.
\end{proof}

\subsection{Nonuniqueness modulo circle isometries}
\label{subsec:nonisometry567}

\begin{proof}[Proof of Theorem~\ref{thm:nonunique-intro}]
The traded pair is constructed in Propositions~\ref{prop:x7-trade326}
and~\ref{prop:x5tilde-orientation}, and its optimality follows from
Corollaries~\ref{cor:NQ5} and~\ref{cor:NQ67}.  It remains to show that
the corresponding designs in the two pairs are not isometric.

For the $5$-designs, the exact integers in
\eqref{eq:sixth-certificate} and~\eqref{eq:x5tilde-p6} satisfy
\begin{equation}\label{eq:x5-p6-mod13}
 \Hfive^6P_6(X_5)\equiv1\pmod{13},
 \qquad
 \Hfive^6P_6(\widetilde X_5)\equiv7\pmod{13}.
\end{equation}
Squaring the residues gives $1$ and $10$.  The two integers therefore
cannot agree up to sign, and
\[
 |P_6(X_5)|\neq|P_6(\widetilde X_5)|.
\]
For the $7$-designs, the original expanded certificate gives
\begin{equation}\label{eq:x7-p8-mod13}
 \Hfive^8P_8(X_7)\equiv12\pmod{13},
\end{equation}
whereas the traded set has residue $7$ by \eqref{eq:x7tilde-p8-mod13}.
Their squares are again $1$ and $10$, so
\[
 |P_8(X_7)|\neq|P_8(\widetilde X_7)|.
\]
Every circle isometry has the form $z\mapsto\eta z$ or
$z\mapsto\eta\overline z$, with $|\eta|=1$.  It multiplies $P_k$ by
$\eta^k$, possibly after conjugating it, and so preserves $|P_k|$.
The two inequalities rule out an isometry in either case.
\end{proof}

\section{Constructions of optimal rational 4-designs}\label{sec:optimal2304}

For the $2304$-point construction, we use a smaller shell and delete
conjugate pairs according to their signed abscissas.  An exact change of
some deletion signs gives a second $4$-design.  The two sets have
different absolute fifth moments and are therefore nonisometric;
the fourth-moment divisor proves both optimal.  At the end of the
section, disjoint rational rotations give exact-strength-$4$ designs
with $2304r$ points.

\subsection{The smaller Gaussian shell}\label{subsec:H4-shell}

The smaller construction keeps conjugation symmetry but not antipodality.
We reuse the notation for Gaussian factors and exponent vectors, now with
\[
 (\mathfrak p_1,\ldots,\mathfrak p_6)
 =(2+\ii,3+2\ii,4+\ii,5+2\ii,6+\ii,5+4\ii),
 \qquad m_1=\cdots=m_6=1,
\]
and
\begin{equation}\label{eq:H4}
 \Hfour=\prod_{j=1}^6\Norm(\mathfrak p_j)^{m_j}
 =5\cdot13\cdot17\cdot29\cdot37\cdot41
 =48\,612\,265.
\end{equation}
Let
\[
 \Efour=\{-1,0,1\}^6.
\]
For $e=(e_1,\ldots,e_6)\in\Efour$, define
\[
 \gamma_e
 =\prod_{j=1}^6
   \mathfrak p_j^{m_j+e_j}\overline{\mathfrak p_j}^{\,m_j-e_j}
 =A_e+\ii B_e.
\]
Each $\gamma_e$ has $\Norm(\gamma_e)=\Hfour^2$.  The Gaussian primes
$\mathfrak p_j$ and $\overline{\mathfrak p_j}$ have six distinct
rational-prime norms, and no $\mathfrak p_j$ is associated to its
conjugate.  Unique factorization in $\Z[\ii]$ therefore gives
$4\cdot3^6=2916$ distinct unit multiples of the $\gamma_e$.

The axial points are exactly $\{\pm\Hfour,\pm\ii\Hfour\}$, coming
from $e=0$; unique factorization rules out any others.  There are no
diagonal points, since $|A_e|=|B_e|$ would give $\Hfour^2=2A_e^2$
with $\Hfour$ odd.  The remaining $2912$ points form $364$ eight-point
packets.  Order their positive canonical representatives
lexicographically, and write
\[
 \gamma_i=A_i+\ii B_i,
 \qquad 0<A_i<B_i,
 \qquad A_i^2+B_i^2=\Hfour^2,
 \qquad 0\leq i\leq363.
\]
Set
\[
 \Delta_i=A_i^2-B_i^2,
 \qquad w_i=\Delta_i^2,
 \qquad \mathcal W=\sum_{i=0}^{363}w_i.
\]

\subsection{The first optimal four-design}\label{subsec:first-x4}

An antipodal $4$-design is a $5$-design, so
Theorem~\ref{thm:fifth-divisibility} rules out antipodality for a
$2304$-point rational $4$-design.  Conjugation symmetry is still
available.  It makes every Fourier moment real and leaves four real
moment equations to solve.

The four conjugate pairs in $\Orb(\gamma_i/\Hfour)$ are indexed by
the signed integer abscissas $\{\pm A_i,\pm B_i\}$.  Choose a set
$F\subseteq\{0,\ldots,363\}$ of packets to keep in full.  From every
packet with $i\notin F$, remove one conjugate pair.  Its absolute
abscissa and sign are specified by $\tau_i,\sigma_i\in\{\pm1\}$:
set
\[
 \lambda_i=
 \begin{cases}
 A_i,&\tau_i=1,\\
 B_i,&\tau_i=-1,
 \end{cases}
 \qquad
 \lambda_i^2=\frac{\Hfour^2+\tau_i\Delta_i}{2}.
\]
The pair to be removed has signed abscissa $\sigma_i\lambda_i$.
A full packet contributes eight points and every other packet contributes
six.  To obtain the desired cardinality, we need
\[
 8|F|+6(364-|F|)=2304,
 \qquad\text{or equivalently}\qquad |F|=60.
\]
Thus each of the other $304$ packets loses two points.  The eight-or-six
rule is a restriction of our construction, not a necessary property of
rational $4$-designs.

Assume $|F|=60$.  Let $\mathcal T$ contain all four signed abscissas
from each full packet and all but $\sigma_i\lambda_i$ from each other
packet.  These abscissas are distinct: their absolute values determine
the positive complementary ordinate and hence the packet.  It follows
that $|\mathcal T|=4\cdot60+3\cdot304=1152$.
For $t\in\mathcal T$, let $u_t>0$ be the complementary coordinate,
so that $t^2+u_t^2=\Hfour^2$, and put
\begin{equation}\label{eq:x4-deletion-family}
 X(F,\tau,\sigma)=
 \left\{\frac{t+\ii u_t}{\Hfour},\frac{t-\ii u_t}{\Hfour}:
 t\in\mathcal T\right\}.
\end{equation}
For any such selectors, this is a conjugation-invariant set of $2304$
distinct rational circle points.  The moment conditions remain to be
imposed.

Define the four integer residuals
\begin{equation}\label{eq:x4-selector-residuals}
\begin{aligned}
 R_4&=\sum_{i\in F}w_i+3\mathcal W-576\Hfour^4,
 &R_2&=\sum_{i\notin F}\tau_i\Delta_i,\\
 R_1&=\sum_{i\notin F}\sigma_i\lambda_i,
 &R_3&=\sum_{i\notin F}\sigma_i\tau_i\Delta_i\lambda_i.
\end{aligned}
\end{equation}
For $X=X(F,\tau,\sigma)$, the Fourier moments are
\begin{equation}\label{eq:x4-residual-moments}
\begin{aligned}
 P_4(X)&=\frac{4R_4}{\Hfour^4},
 &P_2(X)&=-\frac{2R_2}{\Hfour^2},\\
 P_1(X)&=-\frac{2R_1}{\Hfour},
 &P_3(X)&=\frac{2R_1}{\Hfour}-\frac{4R_3}{\Hfour^3}.
\end{aligned}
\end{equation}
For the fourth-moment identity, use
\[
 \RePart\!\left[\left(\frac{\gamma_i}{\Hfour}\right)^4\right]
 =\frac{2(A_i^2-B_i^2)^2-(A_i^2+B_i^2)^2}{\Hfour^4}
 =\frac{2w_i}{\Hfour^4}-1.
\]
This real part is unchanged by negation, conjugation, or a quarter-turn.
Thus all four conjugate pairs in a packet have the same fourth moment,
independently of $\tau_i$ and $\sigma_i$.  Summing gives
\[
\begin{aligned}
 P_4(X)
 &=8\sum_{i\in F}\left(\frac{2w_i}{\Hfour^4}-1\right)
   +6\sum_{i\notin F}\left(\frac{2w_i}{\Hfour^4}-1\right)\\
 &=\frac{4}{\Hfour^4}
   \left(\sum_{i\in F}w_i+3\mathcal W-576\Hfour^4\right).
\end{aligned}
\]
This explains the choice of $w_i$ as the packet weight.  In angular
coordinates, $\Delta_i/\Hfour^2=\cos(2\theta_i)$ for
$\gamma_i/\Hfour=e^{\ii\theta_i}$.

A full packet has zero first, second, and third moments.  The pair
deleted from packet $i$ contributes
\[
 \frac{2\sigma_i\lambda_i}{\Hfour},\qquad
 \frac{2\tau_i\Delta_i}{\Hfour^2},\qquad
 \frac{4\sigma_i\tau_i\Delta_i\lambda_i}{\Hfour^3}
   -\frac{2\sigma_i\lambda_i}{\Hfour}.
\]
in these three degrees, respectively.  The third expression uses
$\lambda_i^3=(\Hfour^2\lambda_i+\tau_i\Delta_i\lambda_i)/2$.
Subtract these contributions from the moments of the $364$ full packets
to obtain the remaining identities in \eqref{eq:x4-residual-moments}.

The Fourier criterion now says that $X(F,\tau,\sigma)$ is a spherical
$4$-design exactly when $R_4,R_2,R_1,R_3$ vanish.  These are equations
for sums over all packets; no packet is required to satisfy a separate
moment condition.

We solve the equations in their order of dependence.  The choice of $F$
fixes $P_4$, independently of $\tau$ and $\sigma$.  We then choose
$\tau$ to make $P_2$ zero.  Finally, $\sigma$ must satisfy both
$R_1=0$ and $R_3=0$: the first balances the deleted abscissas, and the
second cancels the cubic moment after its linear part has vanished.
A solution at one stage need not extend to the next.

Algorithm~\ref{alg:first-four} searches for $F$ in overlapping
eligibility lists and checks disjointness of the chosen supports.
Two $16$-block searches then find $\tau$ and $\sigma$.
For $\sigma$, a radix encoding combines the two equations into one
exact half-sum.  The final sorted-list join uses the classical
meet-in-the-middle subset-sum principle of Horowitz and Sahni
\cite{HS74}.  The first four-design search uses reservoir sampling for
capped retention, as described by Vitter~\cite{Vitter85}; this should
be distinguished from the first-qualifying-record retention used in
the second sign search.  At each stage, the first successful output
is kept without backtracking.  Exhausting the trial budget is therefore not a
nonexistence result.  The final audit checks the moments and the
identity of Proposition~\ref{prop:x4-fifth-selector}.
The selectors used below define $X_4$; the searches recorded in the
appendix found a different baseline.

\begin{proposition}\label{prop:x4-selector}
There are sets $F,E,G\subseteq\{0,\ldots,363\}$ with $|F|=60$ and
$E,G\subseteq\{0,\ldots,363\}\setminus F$ such that, for $i\notin F$,
\[
 \tau_i=
 \begin{cases}
  1,&i\in E,\\
 -1,&i\notin E,
 \end{cases}
 \qquad
 \lambda_i=
 \begin{cases}
  A_i,&\tau_i=1,\\
  B_i,&\tau_i=-1,
 \end{cases}
 \qquad
 \sigma_i=
 \begin{cases}
  1,&i\in G,\\
 -1,&i\notin G,
 \end{cases}
\]
and
\begin{align}
 \sum_{i\in F} w_i&=576\Hfour^4-3\mathcal W,\label{eq:x4-selector-w}\\
 \sum_{i\notin F}\tau_i\Delta_i&=0,\label{eq:x4-selector-tau}\\
 \sum_{i\notin F}\sigma_i\lambda_i&=0,\label{eq:x4-selector-linear}\\
 \sum_{i\notin F}\sigma_i\tau_i\Delta_i\lambda_i&=0.
 \label{eq:x4-selector-cubic}
\end{align}
\end{proposition}

\begin{proof}
The sets $F,E,G$ are listed in Appendix~\ref{app:certificate-four}.
Substitution into the four sums verifies the identities by integer
arithmetic.  Both sides of \eqref{eq:x4-selector-w} equal
\[
 169\,434\,476\,302\,395\,680\,888\,170\,854\,924\,924.
\]
\end{proof}

Use the selectors of Proposition~\ref{prop:x4-selector} in
\eqref{eq:x4-deletion-family} to define $X_4=X(F,\tau,\sigma)$.
All four residuals vanish.  Equation~\eqref{eq:x4-residual-moments}
and Lemma~\ref{lem:fourier} therefore show that $X_4$ is a $4$-design.
To compute the next moment, write
\[
 p_r=\sum_{t\in\mathcal T}t^r.
\]
Summing the conjugate-pair identities for degrees $1$ through $4$ gives
\[
\begin{aligned}
 P_1(X_4)&=\frac{2p_1}{\Hfour},
 &P_2(X_4)&=\frac{4p_2}{\Hfour^2}-2304,\\
 P_3(X_4)&=\frac{8p_3}{\Hfour^3}-\frac{6p_1}{\Hfour},
 &P_4(X_4)&=\frac{16p_4}{\Hfour^4}-\frac{16p_2}{\Hfour^2}+2304.
\end{aligned}
\]
Their vanishing is equivalent to
\begin{equation}\label{eq:x4-power-sums}
 (p_1,p_2,p_3,p_4)=(0,576\Hfour^2,0,432\Hfour^4).
\end{equation}
The normalized values $p_2/(1152\Hfour^2)=1/2$ and
$p_4/(1152\Hfour^4)=3/8$ are the second and fourth abscissa moments
of uniform measure on the circle.

Let $\mathcal N=\{0,\ldots,363\}\setminus F$.
For $i\in\mathcal N$, the deleted coordinate is $\sigma_i\lambda_i$
and $\lambda_i^2=(\Hfour^2+\tau_i\Delta_i)/2$.

\begin{proposition}
\label{prop:x4-fifth-selector}
For the verified selector data,
\begin{equation}\label{eq:W5}
 \begin{aligned}
 W_5&:=\sum_{i\in\mathcal N}\sigma_i\Delta_i^2\lambda_i\\
 &=-154\,714\,837\,946\,793\,873\,141\,735\,703\,233\,907\,917\,120
 \neq0.
 \end{aligned}
\end{equation}
Moreover,
\begin{equation}\label{eq:x4-fifth-witness}
 \begin{aligned}
 \Hfour^5P_5(X_4)&=-8W_5\\
 &=1\,237\,718\,703\,574\,350\,985\,133\,885\,625\,871\,263\,336\,960
 \neq0.
 \end{aligned}
\end{equation}
Hence the constructed set has exact strength $4$.
\end{proposition}

\begin{proof}
Every odd power sum over a full packet vanishes.  Subtracting the
omitted coordinates gives
\[
 p_5=-\sum_{i\in\mathcal N}\sigma_i\lambda_i^5.
\]
Since
\[
 \lambda_i^5=\frac{\lambda_i}{4}
 \bigl(\Hfour^4+2\Hfour^2\tau_i\Delta_i+\Delta_i^2\bigr),
\]
the two odd selector equations, \eqref{eq:x4-selector-linear} and
\eqref{eq:x4-selector-cubic}, reduce this to $p_5=-W_5/4$.
For $z=(t+\ii u_t)/\Hfour$, expansion gives
\[
 z^5+\overline z^{\,5}
 =\frac{2(16t^5-20\Hfour^2t^3+5\Hfour^4t)}{\Hfour^5}.
\]
Summing and using $p_1=p_3=0$ yields
$\Hfour^5P_5(X_4)=32p_5=-8W_5$.
The selector data give the values in \eqref{eq:W5} and
\eqref{eq:x4-fifth-witness}.
\end{proof}

\begin{remark}
Exact strength $4$ also follows from Theorem~\ref{thm:fifth-divisibility},
which rules out a $2304$-point rational $5$-design on every shell.
We need the explicit fifth moment below to distinguish the two optimal
$4$-designs up to isometry.
\end{remark}

The construction and Proposition~\ref{prop:x4-selector} give $2304$
distinct rational points with $P_k(X_4)=0$ for $1\leq k\leq4$.
Together with $P_5(X_4)\neq0$ from
Proposition~\ref{prop:x4-fifth-selector}, this proves
Theorem~\ref{thm:optimal4}.

The abscissas give the fourth-moment Diophantine formulation discussed
by Cui, Xia, and Xiang~\cite[Section~1.3]{CXX}.  The $2304$ entries below
are indexed by circle points and repeat across conjugate pairs.  The
$1152$ distinct abscissas $t/\Hfour$, $t\in\mathcal T$, satisfy the
same normalized fourth-moment identity, by
\eqref{eq:x4-power-sums}.

\begin{corollary}\label{cor:diophantine}
There exist rational numbers $x_1,\ldots,x_{2304}$ such that every
$1-x_i^2$ is a square in $\Q$ and
\[
 \frac1{2304}\sum_{i=1}^{2304}x_i^4=\frac38.
\]
For the construction in Theorem~\ref{thm:optimal4}, one may write
$x_i=a_i/\Hfour$ with $a_i,b_i\in\Z$ satisfying
\[
 a_i^2+b_i^2=\Hfour^2,
 \qquad
 \sum_{i=1}^{2304}a_i^4=864\Hfour^4,
 \qquad \Hfour=48\,612\,265.
\]
\end{corollary}

\begin{proof}
Write $X_4=\{(x_i,y_i):1\leq i\leq2304\}$.  Since
$x_i,y_i\in\Q$ and $x_i^2+y_i^2=1$, each $1-x_i^2=y_i^2$ is a rational
square.  Apply the design identity to $f(x,y)=x^4$:
\[
 \frac1{2304}\sum_{i=1}^{2304}x_i^4
 =\frac1{2\pi}\int_0^{2\pi}\cos^4\theta\,d\theta
 =\frac38.
\]
Because $X_4\subset \Hfour^{-1}\Z^2$, we can write
$x_i=a_i/\Hfour$ and $y_i=b_i/\Hfour$ with integral numerators.
Clearing denominators gives the claimed identity.
\end{proof}

\subsection{A second optimal design from an exact sigma trade}
\label{subsec:second-x4}

Keep $F$ and $\tau_i$ fixed, and reverse the deletion signs on a subset
of $\mathcal N$.  The design equations remain valid exactly when the
first and third power sums over that subset vanish.
Algorithm~\ref{alg:second-four} searches for a new full sign vector,
compares it with the original to find the flip set, and tests
inequivalence through the absolute fifth moment.  The recorded run
flips $155$ packets of the rediscovered baseline; the certificate for
$X_4$ below flips $157$.

\begin{proposition}
\label{prop:x4-sigma-trade}
There is a set $\mathcal J\subset\mathcal N$ with $|\mathcal J|=157$ such
that
\begin{equation}\label{eq:x4-sigma-trade}
 \sum_{i\in\mathcal J}\sigma_i\lambda_i=0,
 \qquad
 \sum_{i\in\mathcal J}\sigma_i\lambda_i^3=0.
\end{equation}
Define
\[
 \widetilde\sigma_i=
 \begin{cases}
  -\sigma_i,&i\in\mathcal J,\\
  \sigma_i,&i\notin\mathcal J,
 \end{cases}
 \qquad i\in\mathcal N.
\]
For $i\in F$, retain all four signed abscissas $\{\pm A_i,\pm B_i\}$; for
$i\in\mathcal N$, retain the same four integers except
$\widetilde\sigma_i\lambda_i$.  Let $\widetilde{\mathcal T}$ be the resulting set and put
\begin{equation}\label{eq:X4tilde}
 \widetilde X_4=
 \left\{
  \frac{t+\ii u_t}{\Hfour},\frac{t-\ii u_t}{\Hfour}:t\in\widetilde{\mathcal T}
 \right\},
\end{equation}
where $u_t>0$ and $t^2+u_t^2=\Hfour^2$.  Then $\widetilde X_4$ consists of $2304$
distinct rational points and is a spherical $4$-design of exact strength $4$.
Its direct fifth-moment certificate is
\begin{equation}\label{eq:x4tilde-fifth}
 \begin{aligned}
 \widetilde W_5
  &:=\sum_{i\in\mathcal N}\widetilde\sigma_i\Delta_i^2\lambda_i\\
  &=-291\,469\,255\,409\,768\,860\,218\,819\,404\,431\,702\,445\,760,
 \end{aligned}
\end{equation}
and hence
\begin{equation}\label{eq:x4tilde-p5}
 \begin{aligned}
 \Hfour^5P_5(\widetilde X_4)
  &=-8\widetilde W_5\\
  &=2\,331\,754\,043\,278\,150\,881\,750\,555\,235\,453\,619\,566\,080
  \neq0.
 \end{aligned}
\end{equation}
Moreover,
\[
 |X_4\cap\widetilde X_4|=1990,
 \qquad
 |X_4\mathbin{\triangle}\widetilde X_4|=628,
\]
and no rotation or reflection of $S^1$ maps $X_4$ onto $\widetilde X_4$.
Thus $\widetilde X_4$ is a second globally cardinality-optimal rational
$4$-design, inequivalent to $X_4$ modulo $O(2)$.
\end{proposition}

\begin{proof}
The data in Appendix~\ref{app:certificate-second-four} leave $F$ and
$E$ unchanged from $X_4$.  Summing over the listed flip set
$\mathcal J$ verifies \eqref{eq:x4-sigma-trade}.

The baseline selector identities imply
\[
 2\sum_{i\in\mathcal N}\sigma_i\lambda_i^3
 =\Hfour^2\sum_{i\in\mathcal N}\sigma_i\lambda_i
  +\sum_{i\in\mathcal N}\sigma_i\tau_i\Delta_i\lambda_i=0.
\]
Reversing the signs on $\mathcal J$ and applying
\eqref{eq:x4-sigma-trade} gives
\[
 \sum_{i\in\mathcal N}\widetilde\sigma_i\lambda_i=0,
 \qquad
 \sum_{i\in\mathcal N}\widetilde\sigma_i\lambda_i^3=0.
\]
Using $2\lambda_i^2=\Hfour^2+\tau_i\Delta_i$, the second identity yields
\[
 \sum_{i\in\mathcal N}\widetilde\sigma_i\tau_i\Delta_i\lambda_i
 =2\sum_{i\in\mathcal N}\widetilde\sigma_i\lambda_i^3
  -\Hfour^2\sum_{i\in\mathcal N}\widetilde\sigma_i\lambda_i=0.
\]
Thus the odd selector equations are preserved.  The even equations
\eqref{eq:x4-selector-w} and \eqref{eq:x4-selector-tau} are unchanged
because $F$ and $\tau_i$ are fixed.  Repeating the calculation of
\eqref{eq:x4-power-sums} for $\widetilde{\mathcal T}$ gives
\[
 \left(\sum_{t\in\widetilde{\mathcal T}}t,
       \sum_{t\in\widetilde{\mathcal T}}t^2,
       \sum_{t\in\widetilde{\mathcal T}}t^3,
       \sum_{t\in\widetilde{\mathcal T}}t^4\right)
 =(0,576\Hfour^2,0,432\Hfour^4).
\]
The distinctness argument is the same as for $X_4$.  These power-sum
identities show that \eqref{eq:X4tilde} is a $2304$-point rational
spherical $4$-design.

The reconstructed selector gives $\widetilde W_5$ as in
\eqref{eq:x4tilde-fifth}.  Replace $\sigma_i$ by $\widetilde\sigma_i$
in the proof of Proposition~\ref{prop:x4-fifth-selector} to obtain
$\Hfour^5P_5(\widetilde X_4)=-8\widetilde W_5$.
Equation~\eqref{eq:x4tilde-p5} follows, and its nonzero value proves
exact strength $4$.  Optimality is Corollary~\ref{cor:NQ}.

Counting the points in the expanded coordinate sets gives the stated
intersection and symmetric difference.  For nonisometry, compare
\[
 |\Hfour^5P_5(X_4)|
 =1\,237\,718\,703\,574\,350\,985\,133\,885\,625\,871\,263\,336\,960
 \neq |\Hfour^5P_5(\widetilde X_4)|.
\]
A circle isometry has the form $z\mapsto\eta z$ or
$z\mapsto\eta\overline z$, with $|\eta|=1$, and preserves $|P_5|$.
Thus the two designs are not $O(2)$-equivalent.
\end{proof}

Proposition~\ref{prop:x4-sigma-trade} proves
Theorem~\ref{thm:second4-intro}.

\begin{proof}[Proof of Corollary~\ref{cor:NQ}]
Theorem~\ref{thm:universal-divisibility} gives $\NQ(4)\geq2304$, while
Theorem~\ref{thm:optimal4} gives $\NQ(4)\leq2304$.
\end{proof}

Disjoint rational rotations of $X_4$ give every positive cardinality
divisible by $2304$, without increasing the exact strength.

\begin{corollary}
\label{cor:x4-family}\label{thm:x4-family}\label{cor:family-exact-strength}
Let
\begin{equation}\label{eq:rho-family}
 \mathfrak p_\ast=7+2\ii,\qquad
 \rho=\frac{\mathfrak p_\ast}{\overline{\mathfrak p_\ast}}
      =\frac{45+28\ii}{53}.
\end{equation}
For every integer $r\geq1$, the union
\begin{equation}\label{eq:x4-family}
 X_4^{(r)}=\bigcup_{j=0}^{r-1}\rho^jX_4
\end{equation}
is disjoint and consists of $2304r$ rational points on $S^1$.
Every $X_4^{(r)}$ has exact strength $4$.
\end{corollary}

\begin{proof}
Let $v_{\mathfrak p_\ast}$ be the normalized $\mathfrak p_\ast$-adic
valuation on $\Q(\ii)^\times$.  The norm $\Norm(\mathfrak p_\ast)=53$
does not occur among the shell primes $5,13,17,29,37,41$ of $X_4$.
Hence $v_{\mathfrak p_\ast}(x)=0$ for every $x\in X_4$, whereas
$v_{\mathfrak p_\ast}(\rho)=1$.  An equality $\rho^jx=\rho^ky$ with
$x,y\in X_4$ therefore forces $j=k$ by taking valuations, and then
$x=y$.  The union is disjoint and has the claimed cardinality.
Rationality follows from $\rho\in S^1\cap\Q(\ii)$.

For $1\leq m\leq4$,
\[
 P_m\bigl(X_4^{(r)}\bigr)
 =P_m(X_4)\sum_{j=0}^{r-1}\rho^{jm}=0,
\]
and Lemma~\ref{lem:fourier} gives the $4$-design property.
For exact strength, use $P_5(X_4)\neq0$ from
Proposition~\ref{prop:x4-fifth-selector} and the identity
\[
 P_5\bigl(X_4^{(r)}\bigr)
 =P_5(X_4)\sum_{j=0}^{r-1}\rho^{5j}.
\]
The valuation $v_{\mathfrak p_\ast}(\rho)=1$ rules out $\rho$ being a
root of unity.  Thus the geometric sum is nonzero for every $r\geq1$,
and the strength is exactly $4$.
\end{proof}

\section{Torus immersions and optimal configurations}
\label{sec:applications}

We apply Gromov's construction~\cite{Gro22,Gro25} to the even moments of
$X_5$, obtaining a flat minimal torus immersion with $2304$ distinct
integral frequencies and Euclidean normal curvature equal to Petrunin's
sharp bound~\cite{Pet24}.  The distinction from a construction using repeated
frequencies concerns the Fourier modes, not the dimension of the source:
both sources are two-dimensional.  We then turn to optimal designs
modulo $O(2)$ and their conjectured infinitude at strengths $4$, $5$,
$6$, and $7$.

\subsection{Gromov's torus construction with distinct frequencies}
\label{subsec:gromov-application}

We use Gromov's even-moment construction
\cite[Sections~2--3]{Gro22}; see also his later account
\cite[Section~13]{Gro25} and Chodosh--Li's formulation
\cite[Remark~11]{CL}.  Applying it to $X_5$ gives a specified set of
distinct common-shell frequencies.  Here $n=2$: a design on $S^1$
produces an immersion of $\T^2$.  This is an application of the
construction recalled in Chodosh--Li's remark, not of their main
$S^2\times S^1$ theorem.

Let $B^N\subset\mathbb R^N$ be the closed unit ball.  For an immersion
$f:\T^2\to B^N$ with Euclidean second fundamental form $A_f$, set
\[
 \mathfrak c(f)=\sup\{|A_f(v,v)|:p\in\T^2,\ v\in T_p\T^2,
 |v|_{f^*g}=1\},
 \qquad
 \mathcal C_N(\T^2)=\inf_f\mathfrak c(f).
\]

\begin{theorem}\label{thm:torus-intro}
There exist $2304$ integral vectors
\[
 u_1,\ldots,u_{2304}\in\Z^2,
 \qquad |u_j|=\Hfive,
 \qquad u_j\neq\pm u_k\quad(j\neq k),
\]
such that
\[
 F(\theta)=\frac1{\sqrt{2304}}
 \bigl(e^{\ii\langle u_1,\theta\rangle},\ldots,
       e^{\ii\langle u_{2304},\theta\rangle}\bigr)
\]
defines a flat minimal immersion
\[
 F:\T^2\longrightarrow S^{4607}\subset B^{4608}
\]
with
\[
 F^*g=\frac{\Hfive^2}{2}(d\theta_1^2+d\theta_2^2),
 \qquad
 |A_F(v,v)|=\sqrt{\frac32}
\]
for every induced-metric unit tangent vector $v$.  Consequently
\[
 \mathcal C_N(\T^2)=\sqrt{\frac32}
 \qquad(N\geq4608).
\]
\end{theorem}

No point of $X_5$ has zero first coordinate.  Choose from each antipodal
pair the representative with positive first coordinate, and write
\[
 X_5=Y\sqcup(-Y),
 \qquad
 Y=\{\widehat u_1,\ldots,\widehat u_m\},
 \qquad m=2304.
\]
Set
\[
 u_j=\Hfive\widehat u_j\in\Z^2.
\]
These vectors are distinct up to sign and satisfy $|u_j|=\Hfive$.

\begin{lemma}\label{lem:torus-moments}
For every $v\in\mathbb R^2$,
\begin{align}
 \frac1m\sum_{j=1}^m\langle u_j,v\rangle^2
 &=\frac{\Hfive^2}{2}|v|^2,\label{eq:second-tensor}\\
 \frac1m\sum_{j=1}^m\langle u_j,v\rangle^4
 &=\frac{3\Hfive^4}{8}|v|^4.\label{eq:fourth-tensor}
\end{align}
\end{lemma}

\begin{proof}
The two integrands are even, so their averages over $Y$ and
$X_5=Y\sqcup(-Y)$ agree.  Since $X_5$ is a $5$-design, these are the
circular averages of $\langle x,v\rangle^2$ and
$\langle x,v\rangle^4$, namely $|v|^2/2$ and $3|v|^4/8$.
Rescaling by $\Hfive$ gives the identities.
\end{proof}

\begin{proof}[Proof of Theorem~\ref{thm:torus-intro}]
Identify $\mathbb R^{2m}$ with $\C^m$ and define
\[
 F(\theta)=\frac1{\sqrt m}
 \bigl(e^{\ii\langle u_1,\theta\rangle},\ldots,
       e^{\ii\langle u_m,\theta\rangle}\bigr),
 \qquad
 \theta\in\T^2=\mathbb R^2/(2\pi\Z)^2.
\]
Integrality of the frequencies makes $F$ well-defined on $\T^2$.
Its norm is $|F|=1$.  Polarizing \eqref{eq:second-tensor} gives
\[
 \langle dF(v),dF(w)\rangle=\frac{\Hfive^2}{2}\langle v,w\rangle,
\]
so
\[
 F^*g=\frac{\Hfive^2}{2}(d\theta_1^2+d\theta_2^2).
\]
In particular, $F$ is an immersion and its induced metric is flat.

Differentiating once more gives
\[
 D_vdF(w)=-\frac1{\sqrt m}
 \bigl(\langle u_j,v\rangle\langle u_j,w\rangle
 e^{\ii\langle u_j,\theta\rangle}\bigr)_{j=1}^m.
\]
For any $a\in\mathbb R^2$, take the real Euclidean inner product in
$\mathbb C^m$:
\[
 \bigl\langle D_vdF(w),dF(a)\bigr\rangle
 =\frac1m\sum_{j=1}^m
 \RePart\!\left(
   \ii\,\langle u_j,v\rangle\langle u_j,w\rangle
        \langle u_j,a\rangle
 \right)=0.
\]
Hence $D_vdF(w)$ is normal, so it equals the Euclidean second
fundamental form $A_F(v,w)$.  By \eqref{eq:fourth-tensor},
\[
 |A_F(v,v)|^2
 =\frac1m\sum_j\langle u_j,v\rangle^4
 =\frac{3\Hfive^4}{8}|v|_0^4.
\]
For a unit vector in the induced metric, $|v|_0^2=2/\Hfive^2$.
Substitution gives
\[
 |A_F(v,v)|^2=\frac32.
\]

Minimality follows from Takahashi's eigenfunction characterization of
minimal immersions in a sphere~\cite[Theorem~3]{Tak66}.  We verify the
criterion directly with our sign convention.  All squared frequency
lengths equal $\Hfive^2$, so, with
$\Delta_g=\operatorname{tr}_g\nabla^2$, we have
\[
 \Delta_gF=-2F.
\]
Thus $\operatorname{tr}_g A_F=\Delta_gF=-2F$.
If $A_F^{S}$ denotes the second fundamental form in the unit sphere,
then
\[
 A_F=A_F^{S}-gF.
\]
Taking the $g$-trace and using $\dim\T^2=2$ gives
\[
 -2F=\operatorname{tr}_gA_F^{S}-2F,
\]
so $\operatorname{tr}_gA_F^{S}=0$ and $F$ is minimal in the sphere.

Petrunin's lower bound for immersions of $\T^n$ into the unit ball is
$\sqrt{3n/(n+2)}$~\cite[Theorem~1.1]{Pet24}.  At $n=2$, it equals
$\sqrt{3/2}$, attained by this immersion in $B^{4608}$.
Composing with the standard linear inclusion preserves the curvature
value in $B^N$ for every $N\geq4608$.
\end{proof}

\begin{remark}
\label{rem:torus-distinct}
Gromov traces the rational even-moment identity to Hilbert's proof of
Waring's theorem; see \cite[Section~2.1.B]{Gro22}.  In the unit-sphere
formulation used in \cite[Section~13]{Gro25}, positive rational weights
are cleared by assigning integer multiplicities to rational directions.
This produces a multiset and does not by itself give distinct nodes.
In the Fourier realization, let
$\chi_u(\theta)=e^{\ii\langle u,\theta\rangle}$ be the character associated
with $u\in\Z^2$.  Equal frequencies $u_j=u_k$ give identical characters
$\chi_{u_j}=\chi_{u_k}$, with the same phase and winding pattern.
If a frequency occurs $r$ times in an equal-amplitude list of length
$M$, its contribution to each moment average has weight $r/M$.

These repetitions can be combined by a fixed unitary change of target
coordinates, which sends the block
\begin{equation}\label{eq:repeated-mode-compression}
 \frac1{\sqrt M}\bigl(\chi_u(\theta),\ldots,\chi_u(\theta)\bigr)
 \quad\longmapsto\quad
 \bigl(\sqrt{r/M}\,\chi_u(\theta),0,\ldots,0\bigr),
\end{equation}
where the block on the left has $r$ entries.  Thus repetition at equal
amplitude is equivalent to a single copy with larger amplitude.
This is an ambient isometry.  Removing the zero coordinates changes
neither the induced metric nor the norm of the Euclidean second
fundamental form, and introduces no new Fourier modes.

In our construction, the $m=2304$ frequencies lie on
$A^2+B^2=\Hfive^2$ and satisfy $u_j\neq\pm u_k$ for $j\neq k$.
They give distinct characters at equal amplitude.  For
$u=(u^{(1)},u^{(2)})$, the entries $u^{(1)}$ and $u^{(2)}$ are the
winding numbers around $S^1$ along the two coordinate loops of the
source.  Different frequencies therefore have different winding
patterns.  The image lies in the Clifford torus
\[
 \T^m_{\mathrm{Cl}}=
 \{(w_1,\ldots,w_m)\in\C^m:|w_j|=m^{-1/2}\text{ for all }j\}.
\]
Although there are $m$ phases, all depend on the same two variables
$\theta_1,\theta_2$.  After rescaling the coordinate circles, the map is
a homomorphism from a compact torus with differential of rank two;
its image is therefore a two-dimensional subtorus.  Distinct frequencies
do not by themselves exclude a finite covering of that image.
\end{remark}

\subsection{The number of optimal rational designs}
\label{subsec:optimal-infinitude}

The constructions above show that optimal rational spherical
designs need not be unique up to orthogonal transformations.
It is natural to ask how many such designs exist.

For $t\in\{4,5,6,7\}$, let
\begin{equation}\label{eq:optimal-moduli}
\mathcal M_t=
\left\{X\subset S^1\cap\Q^2:
\begin{array}{l}
X\text{ is an equal-weight spherical }t\text{-design},\\
|X|=\NQ(t)
\end{array}\right\}\big/\!\sim,
\end{equation}
where $X\sim Y$ if $Y=QX$ for some $Q\in O(2)$.

Theorems~\ref{thm:second4-intro} and
\ref{thm:nonunique-intro} show that
$|\mathcal M_t|\ge 2$, while the rationality
of the nodes implies that $\mathcal M_t$ is
at most countable.

\begin{problem}
\label{prob:optimal-infinitude}
For each $t\in\{4,5,6,7\}$, is $\mathcal M_t$
finite or infinite? If it is finite, determine
its exact cardinality.
\end{problem}

\section*{AI statement}

ChatGPT 5.6 Sol was used as an AI-assisted tool in the development
and debugging of the computational code used in this work. The authors
are responsible for verifying and using the code, interpreting its
computational outputs, and ensuring the correctness of all mathematical
statements and conclusions in the paper.

\clearpage
\appendix

% ---- Begin inlined appendix/proof_search_appendix.tex ----
\section{Practical proof-search algorithms}\label{app:proof-search}
\setcounter{algocf}{0}
\renewcommand{\thealgocf}{A.\arabic{algocf}}
\setcounter{equation}{0}
\renewcommand{\theequation}{A.\arabic{equation}}
\renewcommand{\theHequation}{appendix.A.\arabic{equation}}
\setcounter{table}{0}
\renewcommand{\thetable}{A.\arabic{table}}
\renewcommand{\theHtable}{appendix.A.\arabic{table}}

The five searches below produce $Y_t,\widetilde Y_t$, not the fixed
witnesses $X_t,\widetilde X_t$ of the main proofs.  Their sign and orbit
trades have sizes $155$ and $435$, respectively, rather than $157$ and
$326$.  The programs start from the shell factors and earlier search
outputs, not from the manuscript masks.  The accompanying ancillary code archive and
its final certificates are described in Appendix~\ref{app:certificates}.

\subsection{Exact tests and list merging}\label{app:exact-primitives}
For an integer point list $V$ on the shell of radius $H$, write
$M_k(V)=\sum_{v\in V}v^k=H^kP_k(V/H)$.
The test $\Audit(V,H;n,t)$ requires $n$ distinct entries, norm $H^2$ for
each entry, and $M_1=\cdots=M_t=0$.  Exact strength is checked by
$M_{t+1}\ne0$.

To solve $\sum_{i\in J}v_i=T$, leaf records retain their supports and sums,
with $T$ subtracted in exactly one leaf.  A merge keeps
\begin{equation}\label{eq:app-modular-merge}
 a+b\equiv0\pmod p,
 \qquad (a,b)\longmapsto(a+b)/p.
\end{equation}
After $\ell$ levels, each value is the original partial residual divided
by $p_1\cdots p_\ell$.  Thus an exact root match $a+b=0$ yields the target
sum.  The recovered support is checked again in $\Z$, including
cardinality and disjointness.  Equal sums with different supports remain
distinct records, subject to the cap.  Caps and intermediate residue
conditions make these searches incomplete; stages are committed, without
exhaustive backtracking.

\subsection{The two four-designs}\label{app:first-four}\label{app:second-four}
For the $364$ ordered packets of Section~\ref{sec:optimal2304}, set
\[
 \Delta_i=A_i^2-B_i^2,\quad w_i=\Delta_i^2,\quad
 T_F=576H_4^4-3\sum_iw_i,\quad
 \mathcal N=\{0,\ldots,363\}\setminus F.
\]
The coordinate selector solves
$\sum_{i\in E}\Delta_i=\tfrac12\sum_{i\in\mathcal N}\Delta_i$;
put $\tau_i=2\mathbf1_E(i)-1$ and
$\lambda_i=A_i$ if $\tau_i=1$, $\lambda_i=B_i$ otherwise.
The two sign equations are encoded by
\begin{equation}\label{eq:app-sigma-radix}
 B_\sigma=1+2\sum_{i\in\mathcal N}\lambda_i,
 \qquad c_i=\lambda_i+B_\sigma\tau_i\Delta_i\lambda_i.
\end{equation}
Indeed, $|\sum_i\sigma_i\lambda_i|<B_\sigma$, so
$\sum_i\sigma_ic_i=0$ is equivalent to the two unencoded equations.
The plus support $G$ therefore satisfies
$\sum_{i\in G}c_i=\tfrac12\sum_{i\in\mathcal N}c_i$.

% ---- Begin inlined appendix/algorithm_1.tex ----
\begingroup
\SetAlFnt{\small}\SetAlCapFnt{\small}\SetAlCapNameFnt{\small\bfseries}
\begin{algorithm}[H]
\caption{The first four-design}\label{alg:first-four}
\Input{The $H_4$ shell and the parameters in Table~\ref{tab:run-parameters}.}
\Output{$F,E,G$ and $Y_4$, or $\UNKNOWN$.}
Rebuild the $364$ packets and $T_F$\;
Execute the next three searches in order; return $\UNKNOWN$ immediately if a stage exhausts its budget\;
Assign each index to eight of $16$ eligibility lists; use three-element leaf supports up to one quarter of the cap and sample four-element supports thereafter\;
Merge with disjoint-support checks; accept $F$ only if $|F|=60$ and $\sum_{i\in F}w_i=T_F$\;
Partition $\mathcal N$ into $16$ blocks of $19$; use nine- and ten-element leaf subsets to solve the exact $\Delta$ half-sum for $E$\;
Form $\tau,\lambda,c_i$; use all $2^{19}$ subsets per block, three quotient levels and a sorted disk join to solve the exact $c$ half-sum for $G$\;
Set $\sigma_i=2\mathbf1_G(i)-1$; check both unencoded sign equations\;
Expand $V_4$; require $\Audit(V_4,H_4;2304,4)$ and $M_5(V_4)=-8W_5\ne0$\;
Export the selector and point table; \Return $Y_4=V_4/H_4$\;
\end{algorithm}
\endgroup

% ---- End inlined appendix/algorithm_1.tex ----

% ---- Begin inlined appendix/algorithm_2.tex ----
\begingroup
\SetAlFnt{\small}\SetAlCapFnt{\small}\SetAlCapNameFnt{\small\bfseries}
\begin{algorithm}[H]
\caption{A second four-design by changing the signs}\label{alg:second-four}
\Input{A verified selector $F,E,G_0$ from Algorithm~\ref{alg:first-four}.}
\Output{$G_1,J,\widetilde Y_4$, or $\UNKNOWN$.}
Keep $F,E,\tau,\lambda$ fixed and reuse the radix coefficients $c_i$\;
\For{the $12$ alternate sigma seeds in Table~\ref{tab:run-parameters}}{
 Solve the same $16$-block half-sum, with all leaf subsets and first-record retention\;
 \lIf{no new exact support $G_1\ne G_0$ is found}{\Continue}
 Expand $V_4'$ and compute $M_5(V_4')=-8W_5'$\;
 \lIf{$\Audit(V_4',H_4;2304,4)$, $W_5'\ne0$, or $|W_5'|\ne|W_5|$ fails}{\Continue}
 Set $J=G_0\mathbin{\triangle}G_1$; check $0<|J|<304$ and
 $\sum_{i\in J}\sigma_i\lambda_i=\sum_{i\in J}\sigma_i\lambda_i^3=0$\;
 Export the new selector and point table; \Return $\widetilde Y_4=V_4'/H_4$\;
}
\Return $\UNKNOWN$\;
\end{algorithm}
\endgroup

% ---- End inlined appendix/algorithm_2.tex ----

\subsection{The five- and seven-design pairs}\label{app:first-pair}
\label{app:orientation-backend}\label{app:trade-lists}\label{app:second-pair-orientation}
Keep the raw representatives $\gamma_i=A_i+\ii B_i$ in the lexicographic
order of $\Efive$, and write $w_i=(A_iB_i/12)^2$.
The selection uses $128$ blocks: $98$ of size $14$ and $30$ of size $15$.
Nine entries per block enforce $1152=128\cdot9$ selected orbits.
For an increasing selection $(i_r)$, orientation solves
\begin{equation}\label{eq:app-child-subset}
 d_r=|\Delta_{i_r}|,\qquad
 \sum_{r\in J}d_r=\tfrac12\sum_rd_r,\qquad
 \sigma_r=2\mathbf1_J(r)-1.
\end{equation}
Both searches use six quotient-merge levels and an exact root join.
Orientation partitions the selected ranks into $128$ blocks of nine and
allows all $2^9$ subsets per block.  Its signs always refer to the new
selected-rank order; their positive count is not prescribed.
Choose $\beta_r\in\{\gamma_{i_r},\ii\gamma_{i_r}\}$ with
$\RePart(\beta_r^2)=\sigma_rd_r$ and put
$\eta=\sum_r(\RePart\beta_r)(\ImPart\beta_r)/12$.

% ---- Begin inlined appendix/algorithm_3.tex ----
\begingroup
\SetAlFnt{\small}\SetAlCapFnt{\small}\SetAlCapNameFnt{\small\bfseries}
\begin{algorithm}[H]
\caption{The first five-design and its completion}\label{alg:five}
\Input{The $H_5$ shell and first-run parameters in Table~\ref{tab:run-parameters}.}
\Output{$\mathcal S_0,Y_5,Y_7$ and an orientation, or $\UNKNOWN$.}
Rebuild the $1822$ candidates; shuffle them into $128$ blocks of size $14$ or $15$\;
Enumerate nine-element leaf supports; shift one leaf by $-H_5^4$, merge through six quotient levels, and join exactly at the root\;
\lIf{no support is found within the budget}{\Return $\UNKNOWN$}
Check distinct indices, $|\mathcal S_0|=1152$, and $\sum_{i\in\mathcal S_0}w_i=H_5^4$\;
Order the selected indices and solve~\eqref{eq:app-child-subset} using the shared orientation backend\;
Expand the four-point orbits to $V_5$; accept only if $\Audit(V_5,H_5;4608,5)$, $M_6\ne0$, and $\eta$ is odd\;
\lIf{no orientation is accepted within the budget}{\Return $\UNKNOWN$}
Set $V_7=V_5\sqcup\ii V_5$; require $\Audit(V_7,H_5;9216,7)$ and $M_8\ne0$\;
Export support, signs and points; \Return $(Y_5,Y_7)=(V_5/H_5,V_7/H_5)$\;
\end{algorithm}
\endgroup

% ---- End inlined appendix/algorithm_3.tex ----

% ---- Begin inlined appendix/algorithm_4.tex ----
\begingroup
\SetAlFnt{\small}\SetAlCapFnt{\small}\SetAlCapNameFnt{\small\bfseries}
\begin{algorithm}[H]
\caption{A second parent from an independent selection}\label{alg:second-seven}
\Input{Verified $\mathcal S_0,Y_7$ and the second selection parameters.}
\Output{$\mathcal S_1,\mathcal R,\mathcal A,\widetilde Y_7$, or $\UNKNOWN$.}
Rerun the $128$-block selection independently of $\mathcal S_0$\;
\lIf{no support is found within the budget}{\Return $\UNKNOWN$}
Check distinct indices, $|\mathcal S_1|=1152$, and $\sum_{i\in\mathcal S_1}w_i=H_5^4$\;
Set $\mathcal R=\mathcal S_0\setminus\mathcal S_1$ and $\mathcal A=\mathcal S_1\setminus\mathcal S_0$; check $0<|\mathcal R|=|\mathcal A|$ and equal $w$-sums\;
Expand $V_7'$; require $\Audit(V_7',H_5;9216,7)$, $M_8(V_7')\ne0$, and $|M_8(V_7')|\ne|M_8(H_5Y_7)|$\;
Export the support, trade and points; \Return $\widetilde Y_7=V_7'/H_5$\;
\end{algorithm}
\endgroup

% ---- End inlined appendix/algorithm_4.tex ----

% ---- Begin inlined appendix/algorithm_5.tex ----
\begingroup
\SetAlFnt{\small}\SetAlCapFnt{\small}\SetAlCapNameFnt{\small\bfseries}
\begin{algorithm}[H]
\caption{Orienting the second parent}\label{alg:second-five}
\Input{Verified $\mathcal S_1,\widetilde Y_7$ and the orientation parameters.}
\Output{New selected-rank signs and $\widetilde Y_5$, or $\UNKNOWN$.}
Order $\mathcal S_1$ increasingly and recompute the $d_r$ and half-sum target\;
\For{each configured orientation trial}{
 Shuffle the ranks into $128$ blocks of nine; enumerate all leaf subsets and use six quotient levels with an exact root join\;
 \lIf{no exact plus support is found}{\Continue}
 Recover $\sigma_r$, choose $\beta_r$, and expand the four-point orbits to $V_5'$\;
 \lIf{$\Audit(V_5',H_5;4608,5)$, $M_6\ne0$, or odd $\eta$ fails}{\Continue}
 Accept the orientation and leave the loop\;
}
\lIf{no orientation was accepted}{\Return $\UNKNOWN$}
Independently require $V_5'\sqcup\ii V_5'=H_5\widetilde Y_7$ and $|M_6(V_5')|\ne|M_6(H_5Y_5)|$\;
Export the signs and points; \Return $\widetilde Y_5=V_5'/H_5$\;
\end{algorithm}
\endgroup

% ---- End inlined appendix/algorithm_5.tex ----

The independent moment checks also use
\begin{align}
 M_6&=4\sum_r\sigma_rd_r(H_5^4-2304w_{i_r}),
 \label{eq:app-child-sixth}\\
 M_8^{(1)}-M_8^{(0)}
 &=16\cdot1152^2\left(\sum_{i\in\mathcal A}w_i^2-
                         \sum_{i\in\mathcal R}w_i^2\right).
 \label{eq:app-trade-moment-change}
\end{align}
The latter follows by cancelling the constant and linear terms in the
eighth orbit moments.  Also
$\eta\equiv\sum_{i\in\mathcal S}q_i^2=H_5^4\equiv1\pmod2$.
Comparisons of absolute failed moments certify nonisometry; a different
seed alone does not.

\subsection{Parameters and verified outputs}\label{app:search-provenance}
Table~\ref{tab:run-parameters} records the supplied runs.  A.1 keeps separate
random-generator states through its trials; A.2 uses seeds
$20261001+104729a$; the other backends use base seed $+1000003a$.
These are integer seeds, not dates.  Exact repetition of a trajectory also
depends on the C++ standard-library implementation.

\begin{table}[H]
\centering\small
\caption{Base seeds, record caps, and trial limits.}\label{tab:run-parameters}
\begin{tabular}{@{}lrrr@{}}
\toprule
stage & base seed & cap & trials\\
\midrule
A.1: $F$ & $20260918$ & $4\,000\,000$ & $20$\\
A.1: $\tau$ & $20260919$ & $200\,000$ & $20$\\
A.1: $\sigma$ & $20260922$ & $12\,000\,000$ & $20$\\
A.2: $\sigma$ & $20261001$ & $12\,000\,000$ & $12\times1$\\
A.3: selection & $20261010$ & $350\,000$ & $1$\\
A.3: orientation & $20261011$ & $200\,000$ & $5$\\
A.4: selection & $20261013$ & $350\,000$ & $3$\\
A.5: orientation & $20261011$ & $200\,000$ & $5$\\
\bottomrule
\end{tabular}
\end{table}
The modulus schedules are
\begin{align*}
 F &: (2097143,2097143,2097143),\\
 \tau &: (8191,8191,8191),\\
 \sigma &: (22901,11999989,11999989),\\
 \text{selection} &: (5,350003,350029,350039,350087,350111),\\
 \text{orientation} &: (3,40009,61,67,71,73).
\end{align*}
A.1 uses reservoir retention; A.2 retains the first qualifying records.
The other backends retain the first qualifying records after shuffling
the input lists.  The final sigma joins use sorted temporary files.
Search records use signed $128$-bit or $192$-bit arithmetic; final moment
tests use arbitrary-precision integers.

\begin{table}[H]
\centering\small
\caption{Failed-moment residues of the rediscovered configurations.
Here $M_k=H^kP_k(Y)$; all displayed moments are real.}
\label{tab:rediscovered-results}
\begin{tabular}{@{}lrrrr@{}}
\toprule
set & points & $k$ & $M_k\bmod101$ & $M_k^2\bmod101$\\
\midrule
$Y_4$ & $2304$ & $5$ & $84$ & $87$\\
$\widetilde Y_4$ & $2304$ & $5$ & $5$ & $25$\\
$Y_5$ & $4608$ & $6$ & $27$ & $22$\\
$\widetilde Y_5$ & $4608$ & $6$ & $55$ & $96$\\
$Y_7$ & $9216$ & $8$ & $74$ & $22$\\
$\widetilde Y_7$ & $9216$ & $8$ & $62$ & $6$\\
\bottomrule
\end{tabular}
\end{table}
In each row $M_1=\cdots=M_{k-1}=0$.  Distinct squared residues separate
the two isometry classes at each strength.  The two five-design
orientations have $587$ and $574$ positive signs.  The previously proved
lower bounds establish cardinality optimality; the seven-designs also
attain the standard minimum at strength six.

% ---- End inlined appendix/proof_search_appendix.tex ----

\clearpage

% ---- Begin inlined appendix/certificates_appendix.tex ----
\section{Final configurations and verification}\label{app:certificates}

The files below belong to the ancillary archive
\path{verification_archive.zip} accompanying this preprint.
After extracting it, witness paths are relative to
\path{data/manuscript/}, and verifier names to
\path{verification/}.  Each configuration has one \path{points.csv}
containing integer numerators and a common denominator; a directory
suffix \path{_second} denotes the corresponding tilde-labelled design.

\subsection{The first five- and seven-designs}\label{app:certificate-pair}
The files \path{X5/orbit_mask.hex} and \path{X5/orientation_mask.hex}
specify Proposition~\ref{prop:selection-certificate}: the former uses
zero-based lexicographic candidate indices, the latter increasing selected
ranks.  The same orbit mask defines $X_7$.
The program \path{verify_manuscript_first.py} reconstructs both sets and
checks their moments, parity and disjoint quarter-turn completion.

\subsection{The traded pair}
\label{app:certificate-traded-parent}\label{app:certificate-traded-child}
The support and trade in Proposition~\ref{prop:x7-trade326} are given by
\path{X7_second/orbit_mask_v2.hex} and \path{X7_second/trade_326.json}.
The child uses \path{X5_second/orientation_mask_v2.hex}, indexed by the
new selected ranks.  The program \path{verify_manuscript_pair.py} checks
the trade, both reconstructions, completion and the absolute-moment
invariants of the two pairs.

\subsection{The two four-designs}
\label{app:certificate-four}\label{app:certificate-second-four}
The files \path{X4/full_packets.txt}, \path{X4/epsilon_plus.txt} and
\path{X4/sigma_plus.txt} encode $F,E,G$ in
Proposition~\ref{prop:x4-selector}; \path{epsilon_plus.txt} means $\tau_i=1$.
The second set shares $F,E$ and uses
\path{X4_second/flipped_sigma_packets.txt} and
\path{X4_second/sigma_plus_v2.txt}.
The program \path{verify_manuscript_first.py} checks the selectors, the
sign trade, both coordinate tables and their fifth-moment separation.

\subsection{Search outputs and commands}\label{app:certificate-running}
The rediscovered selectors are the \path{selector.json} files under
\path{data/rediscovered/}; each five-design selector also determines its
seven-design parent.  They are checked by \path{verify_search_x4.py} and
\path{verify_search_pair.py}.  From the archive root, run
\begin{verbatim}
python verify.py
\end{verbatim}
This also runs \path{audit_final_points.py}, which independently checks all
$12$ coordinate tables, the required monomials and complex moments,
exact strengths, four completions, six nonisometry comparisons and the
torus moment tensors.  It uses only the Python standard library, does not
modify the data, and does not start a search.  To run the searches instead,
use
\begin{verbatim}
bash search/rerun_all.sh /path/to/a/fresh/workspace
\end{verbatim}
The script compiles the four C++ backends and runs Algorithms A.1--A.5
with the parameters of Table~\ref{tab:run-parameters}.  It requires Linux,
a C++17 compiler, Boost headers, and additional memory and temporary disk
space.  New outputs go to the specified workspace.

% ---- End inlined appendix/certificates_appendix.tex ----

\end{document}